\documentclass[11pt,dvips,twoside]{article}

\usepackage{pslatex}
\usepackage{fancyhdr}
\usepackage{graphicx}
\usepackage{geometry}

\RequirePackage{amsfonts,amssymb,amsmath,amscd,amsthm}
\RequirePackage{txfonts}
\RequirePackage{graphicx}
\RequirePackage{xcolor}
\RequirePackage{geometry}
\RequirePackage{enumerate}

\def\figurename{Figure} 
\makeatletter
\renewcommand{\fnum@figure}[1]{\figurename~\thefigure.}
\makeatother

\def\tablename{Table} 
\makeatletter
\renewcommand{\fnum@table}[1]{\tablename~\thetable.}
\makeatother

\usepackage{amsmath}
\usepackage{amssymb}
\usepackage{amsfonts}
\usepackage{amsthm,amscd}

\newtheorem{theorem}{Theorem}[section]

\theoremstyle{definition}
\newtheorem{definition}[theorem]{Definition}

\theoremstyle{remark}
\newtheorem{remark}[theorem]{Remark}

\numberwithin{equation}{section}

\begin{document}

\title{\bfseries\scshape{On wellposedness of generalized neural field equations with delay}}
\author{\bfseries\scshape Evgenii Burlakov\thanks{e-mail address: evgenii.burlakov@nmbu.no}\\
Norwegian University of Life Sciences,\\ Department of Mathematical
Sciences and Technology,
P.O. Box 5003,\\
{\AA}s  1432  Norway \\
\\\bfseries\scshape Evgeny Zhukovskiy\thanks{e-mail address: zukovskys@mail.ru}\\
Tambov State University, \\ Department of Mathematics, Physics and
Computer Sciences,
31~Internatsionalnaya st.,\\
Tambov 392000 Russia \\
\\\bfseries\scshape Arcady Ponosov\thanks{e-mail address: arkadi.ponossov@nmbu.no}\\
Norwegian University of Life Sciences,\\ Department of Mathematical
Sciences and Technology,
P.O. Box 5003,\\
{\AA}s  1432  Norway \\
\\\bfseries\scshape John Wyller\thanks{e-mail address: john.wyller@nmbu.no}\\
Norwegian University of Life Sciences,\\ Department of Mathematical
Sciences and Technology,
P.O. Box 5003,\\
{\AA}s  1432  Norway
\\ \\ \\{\rm (Communicated by Editor-in-Chief)}}

\date{}
\maketitle

\thispagestyle{empty} \setcounter{page}{1}

\thispagestyle{fancy} \fancyhead{}

\fancyhead[L]{{\LARGE J}ournal of {\LARGE A}bstract {\LARGE D}ifferential {\LARGE E}quations and {\LARGE A}pplications\\
Volume 2, Number 1, pp. {\thepage--\pageref{lastpage-01} (2011)}}
\fancyhead[R]
\fancyfoot{}
\renewcommand{\headrulewidth}{.0000pt}

\noindent\hrulefill

\noindent {\bf Abstract.} We obtain conditions for existence of
unique global or maximally extended solutions
 to generalized neural field equations. We also study continuous dependence of these
 solutions on the spatiotemporal integration kernel, delay effects, firing rate and prehistory
 functions.

\noindent \hrulefill

\vspace{.3in}


\vspace{.08in} \noindent \textbf{Keywords}: Neural field equations;
wellposedness; integral equations with delay.

\vspace{.3in}

\noindent {\small Received: Month Day, 2010 $\parallel$  Accepted: Month Day, 2010}

\vspace{.2in}

\section{Introduction}

Firing rate models are used in the investigation of the properties
of   strongly interconnected cortical networks. In neural field
models
  the cortical tissue has in addition been modeled as continuous lines or sheets
  of neurons. In such models the spatiotemporally varying neural activity
  is described by a single or several scalar fields, one for each neuron type
   incorporated in the model. These models are formulated in terms of differential,
integro-differential equations and integral equations. The most
well-known and simplest model in that respect is the Amari model
(see e.g.  \cite{Amari})
\begin{eqnarray}
\begin{array}{c}\displaystyle
u_t(t,x)= -u(t,x)+\int\limits_{R}\omega(x-y)f(u(t,y))dy+I(t,x)+h, \\
t\geq0, x\in R.
\end{array}
 \end{eqnarray}
Here the function $u(t,x)$ denotes the activity of a neural element
at time $t$ and position $x$. The connectivity function (spatial
convolution kernel) $\omega(x)$ determines the coupling between the
elements and the non-negative function $f(u)$ gives the firing rate
of a neuron with activity $u$. Neurons at a position $x$ and time
$t$ are said to be active  if $f(u(t,x))>0$. The function $I(t,x)$
and the parameter $h$ represent a variable and a constant external
inputs, respectively.

The literature on the Amari model (1.1) and its extensions is vast.
The key issues in most of the published papers on these models are
existence and stability of coherent structures like localized
stationary solutions (so-called \emph{bumps}) and traveling
fronts/pulses, pattern formation as the outcome of a Turing type of
instability and issues like wellposedness of the actual models. See
e.g. the reviews  \cite{E1998},  \cite{C2005} and  \cite{B2012} (and
the references therein) for more details.


\null

 \verb"This is a draft of the paper containing the main results with the proofs."

\verb"Full-text version is available at "

http://math-res-pub.org/jadea/6/1/wellposedness-generalized-neural-field-equations-delay

\null

\begin{eqnarray}
\begin{array}{c}\displaystyle
Au_{t}(t,x)=-u(t,x)+\int\limits_{R}\mathbb{W}(x-y)f(u(t,y))dy,\quad
t\geq0, x\in R,
\end{array}
 \end{eqnarray}
$$
 A=\left(%
\begin{array}{cc}
  \!\!1 & 0\!\! \\
  \!\!0 &\alpha\!\! \\
\end{array}%
\right), \!\!\quad
\mathbb{W}(x)=\left(%
\begin{array}{cc}
  \!\!\omega_{ee}(x) & -\omega_{ei}(x)\!\! \\
  \!\!\omega_{ie}(x) & -\omega_{ii}(x)\!\! \\
\end{array}%
\right), \!\!
$$
$$
u(t,x)=\left(%
\begin{array}{c}
  \!\!u_e(t,x)\!\! \\
   \!\!u_i(t,x)\!\! \\
\end{array}%
\right),\!\!\quad
f(u(t,x))=\left(%
\begin{array}{c}
  \!\!f_e(u_e(t,x))\!\! \\
  \!\!f_i(u_i(t,x))\!\! \\
\end{array}%
\right).
$$

\pagestyle{fancy} \fancyhead{} \fancyhead[EC]{E. Burlakov, E.
Zhukovskiy, A. Ponosov and J.Wyller} \fancyhead[EL,OR]{\thepage}
\fancyhead[OC]{On Wellposedness of Generalized Neural Field
Equations with Delay} \fancyfoot{}
\renewcommand\headrulewidth{0.5pt}

\begin{eqnarray}
\begin{array}{c}\displaystyle
u_t(t,x)=-Lu(t,x)+
\int\limits_{\Omega}\omega(t,x,y)f(u(t-\tau(x,y),y))dy + I(t,x),\\
 t\in[a,\infty),\ x\in \Omega\subset R^m
\end{array}
 \end{eqnarray}
%

 \begin{eqnarray}
\begin{array}{c}\displaystyle
u^\varepsilon_t(t,x)=
-u^\varepsilon(t,x)+\int\limits_{R}\omega^\varepsilon(x-y)f(u^\varepsilon(t,y))dy,
\\ t\geq0, x\in R,\end{array}
 \end{eqnarray}

 \begin{eqnarray}
\begin{array}{c}\displaystyle
u_t(t,x_c,x_f)=-u_0(t,x_c,x_f)+\int\limits_{R^m}\int\limits_{\mathcal{Y}}
\omega(x_c-y_c,x_f-y_f)f(u(t,y_c,y_f)) dy_c dy_f, \\ t>0,\ x_c\in
R^m,\ x_f\in \mathcal{Y}\subset R^m.\end{array}
 \end{eqnarray}

The Volterra formulation
 \begin{eqnarray}
\begin{array}{c}\displaystyle
u(t,x)=\int\limits_{-\infty}^t\eta(t-s)
\int\limits_{R}\omega(x-y)f(u(s-|x-y|/v,y))dyds, \\ t\in R,\ x\in
R.\end{array}
 \end{eqnarray}
has been investigated by Venkov 
 \begin{eqnarray}
\begin{array}{c}\displaystyle
u(t,x)=\int\limits_{-\infty}^t  \int\limits_{\Omega}W(t,s,x,y)
f(u(s-\tau(s,x,y),y))dyds,\\
 t\in R,\ x\in \Omega \subseteq R^m\end{array}
 \end{eqnarray}
 \begin{eqnarray}
\begin{array}{c}\displaystyle
u(t,x)=\int\limits_a^t \int\limits_{\Omega} W(t,s,x,y)
f(u(s-{\tau}(s,x,y),y))dyds,\\
 t\in[a,\infty), x\in\Omega;   \\
u(\xi,x)\equiv 0,\ \xi\leq a,  x\in\Omega.\end{array}
 \end{eqnarray}

We do not consider  external inputs $I(t,x)$ and $h$ (unlike
 \cite{Amari},  \cite{Faye}) in our models, as they do not involve any
nonlinearities and, hence, only make statements and proofs more
cumbersome. We stress, however, that all the results below remain
valid in the presence of the external inputs as well.

Note that we get $(1.2)$ from $(1.8)$ by taking
$$W(t,s,x,y)=\eta(t,s) \omega(x-y)$$ with
$$\eta(t,s)=diag\Big(\exp\big(-(t-s)\big), \alpha
\exp\big(-\alpha(t-s)\big)\Big) \ \mbox{and} \ \tau(t,x,y)\equiv0.
$$
 If we neglect $I(t,x)$ in $(1.3)$, we can obtain $(1.3)$ from $(1.8)$ with  $$W(t,s,x,y)=\eta(t,s)
\omega(t,x,y),$$ $$\eta(t,s)=diag\Big(l_1 \exp\big(-l_1(t-s)\big),
\ldots, l_n \exp\big(-l_n(t-s)\big)\Big), \ \tau(t,x,y)=\tau(x,y).$$
Taking $\Omega=R^m\times \mathcal{Y}$ ($\mathcal{Y}$ is some
$m$-dimensional torus  \cite{Vis}), $$x=(x_c,x_f), \ y=(y_c,y_f),$$
$$W(t,s,x,y)=\exp\big(-(t-s)\big) \omega(x_c-y_c,x_f-y_f)$$ in $(1.8)$  with
$$\tau(t,x,y)\equiv0,$$ we get the model $(1.5)$. Finally, with
$$W(t,s,x,y)=\eta(t-s) \omega(x-y)$$ and $$\tau(t,x,y)=|x-y|/v$$ in
$(1.7)$, we obtain $(1.6)$, which covers, in turn, the model $(1.1)$
without the external inputs.

 Our results generalize the results
obtained by Potthast \textit{et al} \cite{PG} and Faye \textit{et
al} \cite{Faye} concerning existence of a unique solution to the
Amari model $(1.1)$ in the Banach space of continuous bounded
functions and to the model $(1.3)$ in the space of square integrable
functions on a bounded domain, respectively. Here we also study
dependence of solutions on the parameters.

The paper is organized in the following way. Section 2 is devoted to
the study of
 local solvability, extendability and continuous dependence of solutions to operator
 Volterra equations on parameters. Building on these general results we investigate
 the models $(1.7)$ and $(1.8)$ in Section 3. Section 4 contains conclusions and an
 outlook.

 We stress that one of the challenging parts of out study is application of
 the general theory of Volterra operators to the integral equations $(1.7)$ and $(1.8)$,
 which are defined on
 unbounded spatial and temporal domains. This general setting requires some
 conditions which are difficult to verify (see main theorems in Section~3). In two
 special cases, which are highly relevant for the neural field theory, we can however relax
 these conditions. The analogues of the main theorems for these special cases are
 formulated as remarks in Section 3 and their proofs are given in Appendix.

\section{Existence, uniqueness and continuous dependence of
solutions on parameters: the case of Volterra  operator equations}

Let us introduce the following  notation:

 \begin{description}
    \item[ ]  $R^n$ is the space of vectors consisting of $n$ real
components with the norm  $|\cdot|$;
\item[ ]  $\Omega$ is some closed subset  of $R^m$;
 \item[ ]  $\mathcal{B}$ is some Banach space   with the norm  $\|\cdot\|_\mathcal{B}$;
    \item[ ]  $Y([a,b],
\mathcal{B})$ is a Banach space of functions
   $y:[a,b]  \to \mathcal{B}$ with the norm
  $\|\cdot\|_{Y}$;
  \item[ ]  $\mathfrak{B}(\Omega,
R^n)$ is some Banach space of functions
   $v:\Omega \to R^n$ with the norm
  $\|\cdot\|_{\mathfrak{B}(\Omega,R^n)}$;
    \item[ ] $\Lambda$ is some metric space;
    \item[ ] $\mu$ is the
  Lebesgue measure;
    \item[ ]  $L^p(\Omega, \mu, R^n)$ is the space of all measurable
   and integrable with $p$-th degree
functions $\chi:\Omega\rightarrow R^n$  with the norm
$\|\chi\|_{L^p([\Omega],\mu,R^n)}=\Big(\int\limits_\Omega|\chi(s)|^pds\Big)^{1/p}$,
$1\leq p <\infty$;
 \item[ ]  $BC(\Omega, R^n)$ is the space of all continuous bounded
 functions  $\vartheta:\Omega\rightarrow R^n$  with the
norm $\|\vartheta\|_{BC(\Omega, R^n)}=\sup\limits_{x\in
\Omega}|\vartheta(x)|$;
 \item[ ]  $C_{0}(\Omega,R^n)$ is the space of all
continuous   functions  $\hat{\vartheta}:\Omega\rightarrow R^n$
satisfying the additional condition
$\lim\limits_{|x|\to\infty}|\hat{\vartheta}(x)|=0$
  in the case if $\Omega$ is unbounded, with the norm
$\|\hat{\vartheta}\|_{C_{0}(\Omega,R^n)}=\max\limits_{x\in
\Omega}|\hat{\vartheta}(x)|$;
 \item[ ]  $C([a,b],\mathfrak{B}(\Omega,
R^n))$ is the space of all continuous functions $\nu:[a,b]
\rightarrow \mathfrak{B}(\Omega, R^n)$,
 with the norm
$\|\nu\|_{C([a,b],\mathfrak{B}(\Omega,
R^n))}=\max\limits_{t\in[a,b]}\|\nu(t)\|_{\mathfrak{B}(\Omega,
R^n)}$.
 \item[ ]  $C((-\infty,b],\mathfrak{B}(\Omega,
R^n))$ is the space of all continuous functions
$\hat{\nu}:(-\infty,b] \rightarrow \mathfrak{B}(\Omega, R^n)$ such
that
 $\lim\limits_{t\to-\infty}\|\hat{\nu}(t)\|_{\mathfrak{B}(\Omega,
R^n)}=0$,
 with the norm
$\|\hat{\nu}\|_{C((-\infty,b],\mathfrak{B}(\Omega,
R^n))}=\max\limits_{t\in(-\infty,b]}\|\hat{\nu}(t)\|_{\mathfrak{B}(\Omega,
R^n)}$.
 \end{description}

 In  the notation for
functional spaces we will not indicate the definition domains and
the image sets of functions, provided that this leads to no
ambiguity.

\null

\begin{definition} An operator $\Psi:Y\to Y$ is said to be a
\textit{Volterra operator} (in the sense of A.N. Tikhonov
\cite{Tikh}) if for any $\xi\in(0, b{-}a)$ and any $y_1, y_2\in Y$
the fact that $y_1(t)=y_2(t)$ on $[a,a{+}\xi]$ implies that $(\Psi
y_1)(t)=(\Psi y_2)(t)$ on $[a,a{+}\xi]$. \end{definition}

\null

In what follows we assume that in the space $Y$ the following
condition is fulfilled:

\null

${\mathcal{V}}$\textit{-condition} \cite{Zhuk}: For arbitrary $y\in
Y$, $\{y_i\}\subset Y$ such that $\|y_i-y\|_{Y}\rightarrow0$ and for
any $\xi\in(0, b{-}a)$ if $y_i(t)=0$ on $[a, a{+}\xi]$, then
$y(t,x)=0$ on $[a, a{+}\xi]$.

\null

For any $\xi\in(0, b{-}a)$ let $Y_\xi=Y([a,a{+}\xi], \mathcal{B})$
denote the linear space of restrictions $y_\xi$ of functions $y\in
Y$ to $[a,a{+}\xi]$ which implies that for  each  $y_\xi\in Y_\xi$
there exists at least one extension $y\in Y$ of the function
$y_\xi$.
 Then we can define the norm of $Y_{\xi}$ by $\|y_\xi\|_{Y_{\xi}}= \inf\|y\|_{Y}$,
where the infimum is taken over all extensions $y\in Y$ of the
function $y_\xi$. Hence,  the space $Y_\xi$ becomes a Banach space.

For an arbitrary $\xi\in(0, b{-}a)$ let an operator
$P_\xi:Y\rightarrow Y$ takes each $y_\xi\in Y_{\xi}$ to some
extension $y\in Y$ of $y_\xi$. Moreover, we define the operators
$E_\xi:Y\rightarrow Y_\xi$ by
  $(E_\xi y)(t)=y(t)$,
$t\in[a,a{+}\xi]$ and $\Psi_\xi:Y_\xi\rightarrow Y_\xi$ by $\Psi_\xi
y_\xi=E_\xi \Psi P_\xi y_\xi$, respectively. Note that for any
Volterra operator $\Psi:Y\rightarrow Y$ the operator
$\Psi_\xi:Y_\xi\rightarrow Y_\xi$ is also a Volterra operator and it
is independent of the way $y=P_\xi y_\xi$ extends $y_\xi$.

\null

\textbf{ Definition 2.2.} A Volterra operator $\Psi:Y\rightarrow Y$
is called \textit{locally contracting} if there exists $q < 1$ such
that for any $r > 0$ one can find ${\delta}>0$ such that the
following two conditions are satisfied for all   $y_1, y_2\in Y$,
such that $\|y_1\|_{Y}\leq r$, $\|y_2\|_{Y}\leq r$:
\\

$\mathfrak{q}_1)$  $\ \ \ \ \ \ \ \ \ \ \ \ \ \ \  \ \ \
\|E_{\delta} \Psi y_1-E_{\delta} \Psi y_2\|_{Y_{\delta}}\leq
q\|E_{\delta} y_1-E_{\delta} y_2\|_{Y_{\delta}}$,
\\

$\mathfrak{q}_2)$ for any $\gamma\in(0, b{-}a{-}\delta]$, the
condition $E_\gamma y_1=E_\gamma y_2$ implies that
$$\ \|E_{\gamma{+}{\delta}}\Psi y_1- E_{\gamma{+}{\delta}}\Psi
y_2\|_{Y_{\gamma{+}{\delta}}}\leq q\|E_{\gamma{+}{\delta}}
y_1-E_{\gamma{+}{\delta}} y_1\|_{Y_{\gamma{+}{\delta}}}.$$

\null

The class of locally contracting operators is rather wide. It
includes not only contracting operators, but also, e.g.
$\tau$-Volterra operators.

\null

\textbf{Definition 2.3.} An operator  $\Psi:Y\rightarrow Y$ is
called $\tau$\textit{-Volterra}  if for any $y_1, y_2\in Y$ the
condition $(\Psi y_1)(t)=(\Psi y_2)(t)$ holds true on $[a,a{+}\tau]$
and for any $\xi\in[0, b{-}a{-}\tau]$, if $y_1(t)=y_2(t)$ on
$[a,a{+}\xi]$, then $(\Psi y_1)(t)=(\Psi y_2)(t)$ on
$[a,a{+}\xi{+}\tau]$.

\null

 Notice that  $\tau$-Volterra operators satisfy conditions
$\mathfrak{q}_1)$ and $\mathfrak{q}_2)$ with   $q=0$ and
${\delta}=\tau$, which are independent of a choice of $r$.

Let us now  consider the equation \begin{eqnarray} y(t)=(\Psi y)(t),
\ t\in[a,b],
\end{eqnarray}    where   $\Psi :Y\to Y$ is a Volterra operator.

\null

  \textbf{Definition 2.4.}
     We define a \textit{local solution} to Eq. $(2.1)$ on $[a,a{+}\gamma]$,
 $\gamma\in(0,b{-}a)$
to be a function $y_\gamma\in Y_\gamma$ that satisfies  the equation
$\Psi_\gamma y_\gamma=y_\gamma$ on $[a,a{+}\gamma]$. We define a
\textit{maximally extended solution} to Eq. $(2.1)$  on
$[a,a{+}\zeta)$, $\zeta\in(0,b{-}a]$ to be a function
$y_\zeta:[a,a{+}\zeta)\rightarrow
 \mathcal{B}$, whose restriction $y_\gamma$ to $[a, a{+}\gamma]$
 is a local solution of Eq. $(2.1)$ for any $\gamma<\zeta$ and
$\lim\limits_{\gamma\rightarrow\zeta-0}\|y_\gamma\|_{Y_\gamma}=\infty$.
We define a \textit{global solution} to Eq. $(2.1)$ to be a function
$y\in Y$ that satisfies this equation on the entire interval $[a,
b]$.

\null

 Let us now consider the equation
\begin{eqnarray}
y(t)=(F(y,\lambda))(t), \ t\in[a,b]
 \end{eqnarray}
with a  parameter $\lambda\in \Lambda$, where for each $\lambda\in
\Lambda$  a Volterra operator $F(\cdot,\lambda):Y\to Y$ satisfies
the property: $F(\cdot,\lambda_0)=\Psi$ for some $\lambda_0\in
\Lambda$. Our aim is to formulate conditions for existence and
uniqueness of solutions to Eq. $(2.2)$ on a certain fixed set
$[a,a+\xi]\subset [a,b]$ (We, naturally, also apply Definition 4 to
Eq. $(2.2)$ at each fixed $\lambda \in\Lambda$); and convergence of
these solutions to solution to Eq. $(2.1)$ in the norm of $Y_\xi$ as
$\lambda$ approaches $\lambda_0$. This means, that the problem
$(2.2)$ is wellposed.

\null

 \textbf{Definition 2.5.} For any $\lambda
\in\Lambda_0\subseteq \Lambda$, let the Volterra operator
$F(\cdot,\lambda):Y\to Y$ be given. This family of  operators is
called \textit{uniformly locally contracting} if there exist
$q\geq0$ and $\delta>0$, such that for each $\lambda
\in\Lambda_0\subseteq \Lambda$ the operator
$F(\cdot,\lambda):Y\rightarrow Y$ is locally contracting with the
constants $q$ and $\delta$.

\null

The following theorem represents our main tool to study of the
wellposedness of the models $(1.7)$ and  $(1.8)$. Minding future
applications, we formulate this theorem here in a more general form
than it is needed for the classical neural field theory.

 \null

\textbf{  Theorem 2.1.} \textit{Assume that the following two
conditions are satisfied: } \vspace{5pt}

 1) \textit{There is a  neighborhood $U_0$ of $\lambda_0$
where the operators $F(\cdot,\lambda):Y\to Y$, $\lambda\in U_0$
 are uniformly locally contracting;}

\vspace{5pt}

2) \textit{For arbitrary ${y}\in Y$, the mapping $F:Y\times
\Lambda\to Y$ is continuous at $(y,\lambda_0)$.}

\vspace{10pt}

  \textit{Then for each $\lambda\in U_0$, Eq. $(2.2)$ has a unique global or
maximally extended solution, and each local solution is a
restriction of this solution.}

 \textit{ If  Eq. $(2.2)$ has a global solution $y_0$ at $\lambda=\lambda_0$,
 then for each
 $\lambda$ (sufficiently close to $\lambda_0$) it also
  has a global solution $y=y(\lambda)$, and
 $\|y(\lambda)-y_0\|_{Y}\rightarrow0$ as $\lambda\to\lambda_0$.}

    \textit{If  Eq. $(2.2)$
 has a maximally extended solution $y_{0\zeta}$ defined on
 $[a,a{+}\zeta)$ at $\lambda=\lambda_0$,
 then for any
   $\gamma\in(0, \zeta)$ one can find a neighborhood of
   $\lambda_0$ such that for any
 $\lambda$ in this neighborhood
 Eq. $(2.2)$  has a local solution  $y_{\gamma}=y_{\gamma}(\lambda)$ defined on
  $[a, a{+}\gamma]$
 and
  $\|y_{\gamma}(\lambda)-y_{0\gamma}\|_{Y_{\gamma}}\rightarrow0$
 as $\lambda\to\lambda_0$.}

\begin{proof} Choose a fixed $\lambda\in U_0$.
Let $r>0$, $\xi\in(0, b{-}a)$, 
$y_\xi\in Y_\xi$, $\widehat{y}\in Y$. Let $B_{Y}(\widehat{y},r)$
denote the set of functions $y\in Y$ such that
$\|y-\widehat{y}\|_{Y}< r$ and $Y([a,b], \mathcal{B},y_\xi)$ denote
the set of functions $y\in Y$ such that $E_\xi y=y_\xi$. Put
$B_{Y_{([a,b],y_\xi)}}(\widehat{y},r)=B_{Y}(\widehat{y},r)\bigcap
Y([a,b], \mathcal{B},y_\xi)$.

We  construct the solution in the following way. We set
$r_1=(1-q)^{-1}\|F(0,\lambda)\|_Y+1$ and find all
 $\delta>0$ that satisfy the   condition 1) with $r=r_1$.
For $\delta_1\!=\!\frac{1}{2}\sup\{\delta\}$, we have
$$\|E_{\delta_1} F(y,\lambda)-E_{\delta_1}
F(u,\lambda)\|_{Y_{\delta_1}}\leq q\|E_{\delta_1} y-E_{\delta_1}
u\|_{Y_{\delta_1}}$$ at any $y, u\in B_{Y}(0,r_1)$. Then
$F((B_Y(0,r)),\lambda)\subset B_Y(0,r)$ for $B_Y(0,r)$ with $r\geq
r_1$.  By the Banach fixed point theorem ( \cite{Kolm1}, p. 43) the
mapping $F_{\delta_1}(\cdot,\lambda)$ has a fixed point
$y_{\delta_1}$ in the ball $B_{Y_{\delta_1}}(0,r_1)$. This fixed
point is a local solution to Eq. $(2.2)$. Using the Banach theorem,
one can also prove that for arbitrary $\vartheta_1\in(0,\delta_1)$
and any local solution $\widetilde{y}_{\vartheta_1}$ to Eq. $(2.2)$
defined on $[a, a{+}\vartheta_1]$ it holds that
$\widetilde{y}_{\vartheta_1}(t)=y_{\delta_1}(t)$ at all $t\in[a,
a{+}\vartheta_1]$.

Choose $r_2=(1-q)^{-1}\|F(P_{\delta_1}y_{\delta_1},\lambda)\|_{Y}+1$
and find all possible
  $\delta>0$ that satisfy the  condition 1) with $r=r_2$.
For $\delta_2=\frac{1}{2}\sup\{\delta\}$ at any $y, u\in
B_{Y([a,b],y_{\delta_1})}(P_{\delta_1}y_{\delta_1},r_2)$ we have
$$\|E_{\delta_1{+}\delta_2}
F(y,\lambda)P_{\delta_1}y_{\delta_1}-E_{\delta_1{+}\delta_2}
F(u,\lambda)\|_{Y_{\delta_1{+}\delta_2}}\leq
q\|E_{\delta_1{+}\delta_2} y-E_{\delta_1{+}\delta_2}
u\|_{Y_{\delta_1{+}\delta_2}}.$$ According to the Banach theorem
there exists a fixed point $y_{\delta_1{+}\delta_2}$ of the mapping
$F_{\delta_1{+}\delta_2}(\cdot,\lambda)$ in $B_{Y([a,
a{+}{\delta_1{+}\delta_2}],y_{\delta_1})}(E_{\delta_1{+}\delta_2}P_{\delta_1}y_{\delta_1},r_2)$.
This fixed point is a local solution to Eq. $(2.2)$ defined on $[a,
a{+}{\delta_1{+}\delta_2}]$. It is an extension of the local
solution $y_{\delta_1}$. For any $\vartheta_2\in(0,\delta_2)$ and
any local solution $\widetilde{y}_{\delta_1{+}\vartheta_2}$ to Eq.
$(2.2)$ defined on $[a, a{+}\delta_1{+}\vartheta_2]$, it holds that
$\widetilde{y}_{\delta_1{+}\vartheta_2}(t)=y_{\delta_1{+}\delta_2}(t)$
for all $t\in[a, a{+}\delta_1{+}\vartheta_2]$. Next, let us choose
$r_3=(1-q)^{-1}\|F(P_{\delta_1+\delta_2}y_{\delta_1+\delta_2},\lambda)\|_{Y}+1$,
 find all possible
  $\delta>0$ that satisfy the  condition 1) with $r=r_3$ and repeat
  the procedure, etc.

If the norms of the obtained local solutions are uniformly bounded
 by some  $\mathfrak{M}\in R$, then for
$r=\mathfrak{M}+1$ due to the local contractivity of the operator
$F(\cdot,\lambda):Y\rightarrow Y$ we find ${\delta}$ such that
${\delta_i}\geq\frac{\delta}{2}$ at each of the steps described
above. Therefore, in a finite number of steps we will obtain a
unique global solution to Eq. $(2.2)$. But if such $\mathfrak{M}$
 does not exist, then the number of steps becomes infinite. As a
result, we obtain a unique maximally extended solution to Eq.
$(2.2)$.

We now prove the continuous dependence of solutions on a parameter
$\lambda$. Consider the case when, Eq. $(2.2)$ has global solution
$y_0=y(\lambda_0)\in Y$ at $\lambda=\lambda_0$. Let us find
$\delta>0$ satisfying the condition 1) at $r=\|y_0\|_{Y}+1$,
and any $\lambda\in U_0$. For $k=[\frac{b{-}a}{\delta}]+1$ denote
$\Delta_l=l\delta$, $l=1,2, \ldots , k$. Since the condition 2)
holds true, for any $\varepsilon>0$ one can find $\sigma_1>0$ and a
neighborhood $U_1$ such that for each $\lambda\in U_1$ we have
$$\|F(u,\lambda)-F(y,\lambda_0)\|_{Y}<\frac{(1-q)\varepsilon}{6}$$ for
all $u\in Y$ such that $\|u-y\|_{Y}<\sigma_1$. Assume that
$\sigma_1<\frac{(1-q)\varepsilon}{6}$. Let us find $\sigma_2>0$ and
$U_2$ such that for arbitrary $\lambda\in U_2$ it holds that
$$\|F_{\Delta_{k-1}}(u_{\Delta_{k-1}},\lambda)-
F_{\Delta_{k-1}}(y_{\Delta_{k-1}},\lambda_0)\|_{Y_{\Delta_{k-1}}}
<\frac{(1-q)\sigma_1}{6}$$ for all $u_{\Delta_{k-1}}\in
Y_{\Delta_{k-1}}$,
$\|u_{\Delta_{k-1}}-y_{\Delta_{k-1}}\|_{Y{\Delta_{k-1}}}<\sigma_2$.
Assume that $\sigma_2<\frac{(1-q)\sigma_1}{6}$, $U_2\subseteq U_1$.
There exist $\sigma_3>0$ and $U_3$ such that for any $\lambda\in
U_3$ it holds true that
$$\|F_{\Delta_{k-2}}(u_{\Delta_{k-2}},\lambda)-
F_{\Delta_{k-2}}(y_{\Delta_{k-2}},\lambda_0)\|_{Y{\Delta_{k-2}}}
<\frac{(1-q)\sigma_2}{6}$$ for any $u_{\Delta_{k-2}}\in
Y_{\Delta_{k-2}}$,
$\|u_{\Delta_{k-2}}-y_{\Delta_{k-2}}\|_{Y_{\Delta_{k-2}}}<\sigma_3$;
$\sigma_3<\frac{(1-q)\sigma_2}{6}$, $U_3\subseteq U_2$ etc. We
perform $k$ iterations and at the last step find $\sigma_k$ and
$U_k$, $0<\sigma_k<\frac{(1-q)\sigma_{k-1}}{6}$, $U_k\subseteq
U_{k-1}$.

Let $y_{0\Delta_1}$ denote a local solution to Eq. $(2.2)$ at
$\lambda=\lambda_0$,  that is a fixed point of the operator
$F_{\Delta_1}(\cdot,\lambda_0):Y_{\Delta_1}\!\rightarrow\!
Y_{\Delta_1}$. If
$\|u_{\Delta_{1}}{-}y_{0\Delta_{1}}\|_{\!Y_{\Delta_{1}}}
 <\sigma_k$, then
$$\|F_{\Delta_{1}}(u_{\Delta_{1}},\lambda)-F_{\Delta_{1}}(y_{0\Delta_{1}},\lambda_0)\|_{Y_{\Delta_{1}}}
<\frac{(1-q)\sigma_{k-1}}{6}$$  for all $\lambda\in U_k$. Taking
into account the condition $1)$,  we get for any natural number $m$
that
$$\|F^m_{\Delta_{1}}(y_{0\Delta_{1}},\lambda)-y_{0\Delta_{1}}\|_{Y_{\Delta_{1}}}
\leq\|F^m_{\Delta_{1}}(y_{0\Delta_{1}},\lambda)-
F^{m-1}_{\Delta_{1}}(y_{0\Delta_{1}},\lambda)\|_{Y_{\Delta_{1}}}+
\ldots$$
$$\ldots+\|F_{\Delta_{1}}(y_{0\Delta_{1}},\lambda)-y_{0\Delta_{1}}\|_{Y_{\Delta_{1}}}
\leq(q^{m-1}+\ldots+q+1)\frac{(1-q)\sigma_{k-1}}{6}\leq\frac{\sigma_{k-1}}{6}.$$
Due to the convergence of the approximations
$F^m_{\Delta_{1}}(y_{0\Delta_{1}},\lambda)$ to the fixed point
$y_{\Delta_{1}}=y_{\Delta_{1}}(\lambda)$ of the operator
$F_{\Delta_1}(\cdot,\lambda):Y_{\Delta_{1}}\rightarrow
Y_{\Delta_{1}}$ we obtain that
$\|y_{\Delta_{1}}-y_{0\Delta_{1}}\|_{Y_{\Delta_{1}}}
\leq\frac{\sigma_{k-1}}{6}$ for each $\lambda\in U_k$. Further, let
$y_{0\Delta_2}$ be a local solution to Eq. $(2.2)$ at
$\lambda=\lambda_0$ defined on $[a,a{+}\Delta_{2}]\times R^n$. Then,
for all $\lambda\in U_k\subseteq U_{k-1}$ and any $u_{\Delta_2}\in
B_{Y([a,
a{+}\Delta_2],y_{\Delta_{1}})}(y_{0\Delta_{2}},\sigma_{k-1})$ we get
$$\|F_{\Delta_{2}}(u_{\Delta_{2}},\lambda)-y_{0\Delta_{2}}\|_{Y_{\Delta_{2}}}=\|
F_{\Delta_{2}}(u_{\Delta_{2}},\lambda)-
F_{\Delta_{2}}(y_{0\Delta_{2}},\lambda_0)\|_{Y_{\Delta_{2}}}
<\frac{(1-q)\sigma_{k-2}}{6}.$$ Then
$$\|F_{\Delta_{2}}(u_{\Delta_{2}},\lambda)-u_{\Delta_{2}}\|_{Y_{\Delta_{2}}}
<\sigma_{k-1}+\frac{(1-q)\sigma_{k-2}}{6}<\frac{(1-q)\sigma_{k-2}}{3}.$$
For all $m=1,2,\ldots$ we have
$$\|F^m_{\Delta_{2}}(u_{\Delta_{2}},\lambda)-u_{\Delta_{2}}\|_{Y_{\Delta_{2}}}
\leq\|F^m_{\Delta_{2}}(u_{\Delta_{2}},\lambda)-
F^{m-1}_{\Delta_{2}}(u_{\Delta_{2}},\lambda)\|_{Y_{\Delta_{2}}}+
\ldots$$
 $$\ldots+\|F_{\Delta_{2}}(u_{\Delta_{2}},\lambda)-u_{\Delta_{2}}\|_{Y_{\Delta_{2}}}
\leq(q^{m-1}+\ldots+q+1)\frac{(1-q)\sigma_{k-2}}{3}\leq\frac{\sigma_{k-2}}{3}.$$
Taking into account the convergence of the approximations
$F^m_{\Delta_{2}}(u_{\Delta_{2}},\lambda)$ to
$y_{\Delta_{2}}=y_{\Delta_{2}}(\lambda)$  we obtain
$$\|y_{\Delta_{2}}-y_{0\Delta_{2}}\|_{Y_{\Delta_{2}}}\leq
\|y_{\Delta_{2}}-F^m_{\Delta_2}(u_{\Delta_{2}},\lambda)\|_{Y_{\Delta_{2}}}+$$$$+
\|F^m_{\Delta_2}(u_{\Delta_{2}},\lambda)-u_{\Delta_{2}}\|_{Y_{\Delta_{2}}}+
\|u_{\Delta_{2}}-y_{0\Delta_{2}}\|_{Y_{\Delta_{2}}}\leq
\frac{\sigma_{k-2}}{3}+\sigma_{k-1} \leq\frac{\sigma_{k-2}}{2}.$$
Using the convergence of sequential approximations
$F^m_{\Delta_{3}}(u_{\Delta_{3}},\lambda)$ to a fixed point
$y_{\Delta_{3}}=y_{\Delta_{3}}(\lambda)$ of the operator
$F_{\Delta_3}(\cdot,\lambda):Y_{\Delta_3}\rightarrow Y_{\Delta_3}$
for any  $u_{\Delta_3}\in B_{Y([a,
a{+}\Delta_3],y_{\Delta_{2}})}(y_{0\Delta_{3}},\sigma_{k-2})$ and
each $\lambda\in U_k\subseteq U_{k-1}\subseteq U_{k-2}$,
 we obtain the estimate
$\|y_{\Delta_{3}}-y_{0\Delta_{3}}\|_{Y_{\Delta_3}}
\leq\frac{\sigma_{k-3}}{2}$. We, then, repeat this procedure.
 At the $k$-th step we
  prove in an analogous way that  the inequality
$\|y(\lambda)-y_0\|_{Y}<\varepsilon$ holds true for all $\lambda\in
U_k$. Therefore, $\|y(\lambda)-y_0\|_{Y}\rightarrow0$ as
$\lambda\rightarrow\lambda_0$.

Let now a solution $y_{0\eta}$ to Eq. $(2.2)$ at $\lambda=\lambda_0$
be maximally extended. Fix arbitrary $\gamma\in(0, \eta)$ and let
$y_{0\gamma}$ denote the restriction of the solution $y_{0\eta}$ to
$[a,a{+}\gamma]\times R^n$. For the equation $u_\gamma=F_\gamma
(u_\gamma,\lambda_0)$ the function $y_{0\gamma}\in
Y([a,a{+}\gamma]\times \Omega, {R}^n)$ is a global solution. As is
shown above, for all $\lambda$ from some neighborhood of $\lambda_0$
the equations $u_\gamma=F_{\gamma} (u_\gamma,\lambda)$ have global
 solutions
 $y_\gamma(\lambda)$,
 and  $\|y_{\gamma}(\lambda){-}y_{0\gamma}\|_{Y_\gamma}\rightarrow0$
 as $\lambda\rightarrow\lambda_0$. \end{proof}

The proof of Theorem 1 has several corollaries which are summarized
in the following remarks:

\null

\begin{remark}If the constant $\delta$ in the condition $1)$ of  Theorem 2.1
is independent of $r$, then Eq. $(2.2)$ has a global solution. This
is the case e.g. for $\tau$-Volterra operators.
\end{remark}

\null

 \begin{remark}In case of a priori boundedness of the solution,
 it is possible to extend the solution beyond the point $b$ in the same
 way as it was done in the proof of  Theorem 2.1. This will give a unique
 solution defined on $[a, \infty)$.
\end{remark}

\null

Notice that  the existence of a maximally extended solution to Eq.
$(2.2)$ at $\lambda=\lambda_0$ does not guarantee the existence of
maximally extended solutions to eq $(2.2)$ at $\lambda$ arbitrarily
close to $\lambda_0$. The following example illustrates this fact.

\null

 \textbf{Example 2.1.} Let operators $\Phi(\cdot,\lambda):L^1([0, \pi], \mu,
 {R})\rightarrow L^1([0, \pi], \mu,  {R}), \ \lambda\in[0,\pi]$,
be defined as
$$(\Phi(y,\lambda))(t)=\left\{
 \begin{array}{ccl}
    0, & \textrm{if}  \ t\in [0,\lambda);\\
  \Big(\int\limits^{t-\lambda}_0y(s)ds\Big)^2+1, & \textrm{if}   \ ~t\in [{\lambda},\pi].
 \end{array}
\right.
$$

These operators are Volterra operators and satisfy the condition 1)
of
 Theorem 2.1: For $q=\frac{1}{2}$ and any $r>0$ one can choose
$\delta=\frac{1}{4r}$, and condition $1)$  becomes fulfilled for all
$t\in[0,\pi)$ and any $\lambda\in [0,\pi]$). Condition 2) of the
Theorem 2.1 is also fulfilled. The equation $y(t)=(\Phi(y,0))(t), \
t\in[0,\pi]$ has a unique maximally extended solution
$y(t)\!=\!\frac{1}{\cos^2t}$ defined on $[0, \frac{\pi}{2})$. Now,
since  for any $\lambda\in (0,\pi]$ the operator
$\Phi(\cdot,\lambda)$ is a $\tau$-Volterra operator, the equation
$y(t)=(\Phi(y,\lambda))(t), \ t\in[0,\pi]$ is globally solvable for
each $\lambda\in (0,\pi]$.

\null

When analyzing Theorem 2.1, it is natural to ask the question
whether the maximally extended solutions to Eq. $(2.2)$ are defined
on time intervals with arbitrarily small length. The following two
remarks give answers to that question:

\null

\begin{remark}Let the assumptions of Theorem 2.1 be
fulfilled and let there exist some neighborhood
$\widehat{U}\subseteq U_0$ of $\lambda_0$ such that Eq. $(2.2)$ has
maximally extended solutions $y_{\zeta_\lambda}$ defined on
$[a,a{+}\zeta_\lambda)$ for any $\lambda\in \widehat{U}$. Then
  $\inf\limits_{\forall \lambda\in \widehat{U}} \zeta_\lambda>0$.
Since for all $\lambda\in U_0$ operators $F(\cdot, \lambda)$ are
uniformly locally contracting, we  get  $\inf\limits_{\forall
\lambda\in \widehat{U}} \zeta_\lambda>0$.
\end{remark}

\null

\begin{remark}Let the assumptions of Theorem 2.1 be
fulfilled and let
 for $\lambda=\lambda_0$ and some sequence $\lambda_i\subset U_0$
equation  $(2.2)$ have maximally extended solutions $y_{0\zeta}$ and
$y_{\zeta_i}$ defined on
   $[a,a{+}\zeta)$  and
 $[a,a{+}\zeta_i)$, respectively. Then
  $\beta=\min\{\zeta, \, \inf\limits_{\forall i} \zeta_i\}>0$,
 and either $\beta=\zeta$, or $\beta=\zeta_{i_0}$ at some $i_0$.
\end{remark}

The positivity of $\beta$ follows from   Remark 3. Next, we choose
 arbitrary $\varepsilon>0$ and a sequence $\gamma_j\in(0, \beta)$,
$\gamma_j\rightarrow\beta$, $j\rightarrow\infty$. For each
$\gamma_j\in(0, \beta)$ there exists a finite $\sup\limits_{\forall
i}\|y_{i\gamma_j}\|_{Y_{\gamma_j}}$ otherwise $\beta=\gamma_j$. Let
us associate the number $\gamma_1$ with the corresponding local
solution $y_{i_1\gamma_1}$ to Eq. $(2.2)$ at
$\lambda=\lambda_{i_1}$, where $i_1$ is the least number such that
$\max\{\|y_{0\gamma_1}\|_{Y_{\gamma_1}},\sup\limits_{\forall
i}\|y_{i\gamma_1}\|_{Y_{\gamma_1}}\}-\|y_{i_1\gamma_1}\|_{Y_{\gamma_1}}<\varepsilon$;
we associate the number $\gamma_2$ with the corresponding local
solution $y_{i_2\gamma_2}$ to   Eq. $(2.2)$ at
$\lambda=\lambda_{i_2}$, where $i_2$ is the least number such that
$\max\{\|y_{0\gamma_2}\|_{Y_{\gamma_2}},\sup\limits_{\forall
i}\|y_{i\gamma_2}\|_{Y[a,
a{+}\gamma_2]}\}-\|y_{i_2\gamma_2}\|_{Y_{\gamma_2}}<\varepsilon$
etc. We obtain a subsequence $\{i_j\}$ of numbers of local solutions
$y_{i\gamma_j}$ to Eq. $(2.2)$ such that
$\|y_{i_j\gamma_j}\|_{Y_{\gamma_j}}\rightarrow\infty$ as
$j\rightarrow\infty$. If the subsequence $\{i_j\}$ is bounded, then
one can find a number $i_{j_0}$ such that
$\lim\limits_{\gamma\rightarrow\beta-0}\|y_{i_{j_0}\gamma}\|_{Y_{\gamma}}=\infty$,
i.e. $\zeta_{i_{j_0}}=\beta$. Otherwise, using the fact that
  $\|y_{i_j\gamma}-y_{0\gamma}\|_{Y_{\gamma}}\rightarrow0$ as
$j\rightarrow\infty$ for any $\gamma\in(0, \zeta)$ we obtain
$\lim\limits_{\gamma\rightarrow\beta-0}\|y_{0\gamma}\|_{Y_{\gamma}}=\infty$,
i.e. $\zeta=\beta$.

\medskip

\section{  Existence, uniqueness and continuous dependence of
solutions on parameters: the case of neural field equations}

 In this section we
apply the results obtained in the previous section to a class of
nonlinear integral equations, typical representatives of which can
be found in the neural field theory. For the sake of convenience, we
consider the following generalization of the model $(1.8)$:

 \begin{eqnarray}
\begin{array}{c}\displaystyle
u(t,x)=\varphi(a,x)+\int\limits_a^t \int\limits_{\Omega} W(t,s,x,y)
f(u(s-{\tau}(s,x,y),y))dyds,\\
 t\in[a,\infty), x\in\Omega;  \\
u(\xi,x)=\varphi(\xi,x),\ \xi\leq a,  x\in\Omega.\end{array}
 \end{eqnarray}
under the following assumptions on the functions involved:

\vspace{10pt}

$(A1)$ For any  $b>a$, $(t,x) \in [a, b]\times\Omega$,  the function
$W(t,\cdot,x,\cdot):[a,b]\times \Omega\to R^n$ is measurable.

\vspace{5pt}

$(A2)$ For any  $b>a$, at almost all $(s,y) \in [a, b]\times\Omega$,
the function $W(\cdot,s,\cdot,y):[a,b]\times \Omega\to R^n$ is
uniformly continuous.

\vspace{5pt}

$(A3)$ For any  $b>a$, $t\in [a,b]$, $\int\limits_{\Omega}
\sup\limits_{x\in\Omega}\big|W(t,s,x,y)\big| dy \leq G(s),$
 where $G\in L^1([a,b],\mu ,
R^n)$.

\vspace{5pt}

 $(A4)$ The function $f:R^n\to R^n$ is measurable and for
any $\mathrm{r}>0$ one can find $\mathrm{f}_\mathrm{r}>0$, such that
for all $u\in R^n$, $|u|\leq r$,  it holds true that
$|f(u)|\leq\mathrm{f}_\mathrm{r}$.

 \vspace{5pt}

  $(A5)$
The delay function ${\tau}:R\times \Omega\times \Omega\to[0,\infty)$
is continuous on $R\times \Omega\times \Omega$.

 \vspace{10pt}

$(A6)$ The prehistory function $\varphi$ belongs to $C((-\infty,a],
BC(\Omega, R^n))$.

 \vspace{10pt}

The model $(3.1)$ with $\varphi(\xi,x)\equiv 0$ can be obtained from
$(1.7)$ by taking $W(t,s,x,y)={\eta}(t,s)\omega(x,y)$, where, e.g.
  $${\eta}(t,s)=\left\{
 \begin{array}{ccl}
    \kappa \exp\big(-\kappa (t-s)\big), & \textrm{if}   \ ~t\geq a;\\
 0, & \textrm{if}  \ t<a;
 \end{array}
\right.$$ or $${\eta}(t,s)=\left\{
 \begin{array}{ccl}
    \kappa (t-s) \exp\big(-\kappa (t-s)\big), & \textrm{if}   \ ~t\geq a;\\
 0, & \textrm{if}  \ t<a
 \end{array}
\right.$$ and $\omega$ can be represented by the "Mexican hat"\
 $$\omega(x,y)=M \exp(-m|x-y|)-K \exp(-k|x-y|)$$ or the "wizard hat"\
 $$\omega(x,y)=M(1-|x-y|) \exp(-m|x-y|),$$ and $$f(u)=\left\{
 \begin{array}{ccl}
    u^\kappa/(\theta^\kappa + u^\kappa), & \textrm{if}   \ ~u\geq 0;\\
 0, & \textrm{if}  \ u<0,
 \end{array}
\right.$$ for some  $\kappa>0$, $\theta>0$, $M>K>0$, and $m>k>0$.
These functions satisfy   the conditions $(A1)$ -- $(A4)$. The
condition $(A4)$ is also fulfilled  e.g. for the sigmoidal functions
$$f(u)=\frac{1}{2}\Big(1+\tanh\big(\kappa(u-\theta)\big)\Big)$$ or
$$f(u)=\frac{1}{1+\exp\big(-\kappa(u-\theta)\big)}$$ with some positive
$\kappa$ and $\theta$. We do not assume in $(A4)$ that function $f$
is bounded (as  in the classical neural field theory), because it
allows us to obtain more general results which may have other
applications. If we take the delay functions
${\tau}(t,x,y)=|x-y|/\upsilon$ for some positive velocity $\upsilon$
or ${\tau}(t,x,y)=d(x,y)$ with continuous function $d:R\times R\to
[0, \infty)$ from   \cite{Venkov} and  \cite{Faye}, respectively, we
find out that the condition $(A5)$ is also satisfied.

We introduce the definition of local, maximally extended and global
solutions just as in the previous section (Definition 2.4).

\null

 \textbf{Definition 3.1.}
     We define a \textit{local solution} to Eq. $(3.1)$  on $[a,a{+}\gamma]\times R^n$,
 $\gamma\in(0,\infty)$
to be a function $u_\gamma\in C([a,a{+}\gamma],BC(\Omega, R^n))$
that satisfies the equation $(3.1)$ on $[a,a{+}\gamma]\times R^n$.
We define a \textit{maximally extended solution} to Eq. $(3.1)$ on
$[a,a{+}{\zeta})\times\Omega$, ${\zeta}\in(0,\infty)$ to be a
function $u_{\zeta}:[a,a{+}{\zeta})\times\Omega\rightarrow
 {R}^n$, whose restriction $u_\gamma$ to $[a, a{+}\gamma]\times\Omega$
for any $\gamma<{\zeta}$ is a local solution of Eq. $(3.1)$ and
$\lim\limits_{\gamma\rightarrow{\zeta}-0}\|u_\gamma\|_{C([a,a{+}\gamma],BC(\Omega,
R^n))}=\infty$. We define a \textit{global solution} to Eq. $(3.1)$
to be a function $u:[a,\infty)\times\Omega\to R^n$, whose
restriction $u_\gamma$ to $[a, a{+}\gamma]\times\Omega$  is its
local solution for any $\gamma\in(0,\infty)$.

\null

 \textbf{  Theorem 3.1.}
\textit{Let the assumptions $(A1)$ -- $(A6)$ hold true. If for any
$r>0$ there exists $\widetilde{f}_r\in R$ such that for all $u_1,
u_2 \in R^n$, $|u_1|\leq r$, $|u_2|\leq r$, we have
$|f(u_1)-f(u_2)|\leq \widetilde{f}_r|u_1-u_2|$, then
 Eq.
$(3.1)$ has a unique global or maximally extended solution and each
local solution is a restriction of this solution.}

 \begin{proof}
We will use  Theorem 2.1, namely, the condition 1), which is
responsible for solvability of the Eq. $(2.2)$) and Remark 2.2 of
the previous section to prove the solvability of $(3.1)$.

First,  we choose an arbitrary $b\in(a,\infty)$, define the
following operator
 \begin{eqnarray}
\begin{array}{c}\displaystyle
(Fu)(t,x)=\varphi(a,x)+\int\limits_a^t
\int\limits_{\Omega}W(t,s,x,y)f\Big(\big(S_{\tau}
u\big)(s,x,y)\Big)dyds,\end{array}
 \end{eqnarray}
 $$
(S_{\tau}^\varphi u\big)(t,x,y)=\left\{
 \begin{array}{ccl}
  \varphi(t-{\tau}(t,x,y),x), \!\!\!&\mbox{if}\!\!\!&~t-{\tau}(t,x,y)<a;\\
  u(t-{\tau}(t,x,y),y),\!\!\!&\mbox{if}\!\!\!&~t-{\tau}(t,x,y)\geq a,\\
 \end{array}
\right.
$$
and show that
$$F:C([a,b],BC(\Omega, R^n))\to
C([a,b],BC(\Omega, R^n)).
$$
For any $t\in[a,b]$ and $u\in C([a,b],BC(\Omega, R^n))$ we have
$$
|(Fu)(t,x_1)-(Fu)(t,x_1)|\leq
$$
$$
\leq\Big|\varphi(a,x_1)+\int\limits_a^t
\int\limits_{\Omega}W(t,s,x_1,y)f\Big(\big(S_{\tau}
u\big)(s,x_1,y)\Big)dyds-\ \ \ \ \ \ \ \ \ \ \ \ \ \ \ \ \ \ \ \ \
$$
$$\ \ \ \ \ \ \ \ \ \ \ \ \ \ \ \ \ \ \ \ \
-\varphi(a,x_2)+\int\limits_a^t
\int\limits_{\Omega}W(t,s,x_2,y)f\Big(\big(S_{\tau}
u\big)(s,x_2,y)\Big)dyds\Big|\leq
$$
$$
\leq|\varphi(a,x_1)-\varphi(a,x_2)|+
$$
$$
+\int\limits_a^t
\int\limits_{\Omega}\Big|W(t,s,x_1,y)-W(t,s,x_2,y)\Big|\Big|f\Big(\big(S_{\tau}
u\big)(s,x_1,y)\Big)\Big|dyds+
$$
$$
+\int\limits_a^t
\int\limits_{\Omega}\Big|W(t,s,x_2,y)\Big|\Big|f\Big(\big(S_{\tau}
u\big)(s,x_1,y)\Big)-f\Big(\big(S_{\tau}
u\big)(s,x_2,y)\Big)\Big|dyds.
$$
By the virtue of the assumption $(A6)$, the first term goes to 0 as
$|x_1-x_2|\to0$. The assumptions $(A2)$ -- $(A4)$ and $(A6)$
guarantee convergence to 0 of the second term on the right hand side
of this inequality as $|x_1-x_2|\to0$. The superposition
$f\Big(\big(S_{\tau} u\big)(s,\cdot,y)\Big)$ is continuous as the
assumptions $(A4)$ -- $(A6)$ hold true. This fact and the assumption
$(A3)$ imply convergence of the last term to 0 as $|x_1-x_2|\to0$.
This proves continuity of $(Fu)(t,\cdot)$.

For each $t\in[a,b]$ and any $u\in C([a,b],BC(\Omega, R^n))$ the
function $(Fu)(t,\cdot)$ is bounded by the virtue of the assumptions
(A3), (A4) and (A6).

Finally, we  choose an arbitrary $u\in C([a,b],BC(\Omega, R^n))$
and, assuming that $t_2>t_1$, check that $(Fu)(\cdot,x)$ is
continuous:
$$
\sup\limits_{x\in\Omega}|(Fu)(t_1,x)-(Fu)(t_2,x)|\leq
$$
$$
\leq\sup\limits_{x\in\Omega}\Big|\int\limits_a^{t_1}
\int\limits_{\Omega}W(t_1,s,x,y)f\Big(\big(S_{\tau}
u\big)(s,x,y)\Big)dyds-\ \ \ \ \ \ \ \ \ \ \ \ \ \ \ \ \ \ \ \ \
$$
$$
\ \ \ \ \ \ \ \ \ \ \ \ \ \ \ \ \ \ \ \ \ -\int\limits_a^{t_2}
\int\limits_{\Omega}W(t_2,s,x,y)f\Big(\big(S_{\tau}
u\big)(s,x,y)\Big)dyds\Big|\leq
$$
$$
\leq\sup\limits_{x\in\Omega}\Big|\int\limits_a^{t_1}
\int\limits_{\Omega}\Big(W(t_1,s,x,y)-W(t_2,s,x,y)\Big)f\Big(\big(S_{\tau}
u\big)(s,x,y)\Big)dyds\Big|+
$$
$$
+\sup\limits_{x\in\Omega}\int\limits_{t_1}^{t_2}
\int\limits_{\Omega}\Big|W(t_2,s,x,y)f\Big(\big(S_{\tau}
u\big)(s,x,y)\Big)dyds\Big|\leq
$$
$$
\leq\int\limits_a^{t_1}
\int\limits_{\Omega}\sup\limits_{x\in\Omega}\Big|W(t_1,s,x,y)-W(t_2,s,x,y)\Big|
\sup\limits_{x\in\Omega}\Big|f\Big(\big(S_{\tau}
u\big)(s,x,y)\Big)\Big|dyds+
$$
$$
+\int\limits_{t_1}^{t_2}
\int\limits_{\Omega}\sup\limits_{x\in\Omega}\Big|W(t_2,s,x,y)\Big|\sup\limits_{x\in\Omega}\Big|f\Big(\big(S_{\tau}
u\big)(s,x,y)\Big)\Big|dyds.
$$
We note that by the virtue of the  assumptions $(A2)$ -- $(A4)$ and
$(A6)$, the first term converges to 0 as $t_1-t_2\to0$. The second
summand goes to 0 as the assumptions $(A3)$, $(A4)$ and $(A6)$ hold
true and $t_1-t_2\to0$.

Thus we proved that $F:C([a,b],BC(\Omega, R^n))\to
C([a,b],BC(\Omega, R^n)) $.

Next, we examine the fulfilment  of Theorem 2.1 condition for the
defined above operator $F:C([a,b],BC(\Omega, R^n))\to
C([a,b],BC(\Omega, R^n))$. Choose an arbitrary $q_0\!<\!1$, $r>0$.
Let $\gamma\!\in\!(0, b-a)$ and $u_1(t,\cdot) = u_2(t,\cdot)$,
$t\in[a,a{+}\gamma]$, where $\|u_1\|_{C([a,b],BC(\Omega, R^n))}\leq
r$ and $\|u_2\|_{C([a,b],BC(\Omega, R^n))}\leq r$.
 By assumption, we get the estimates
$$
\sup\limits_{t\in[a,a+\gamma+\delta],
x\in\Omega}\Big|\int\limits_a^{t}
\int\limits_{\Omega}W(t,s,x,y)f\Big(\big(S_{\tau}^\varphi
u_1\big)(s,x,y)\Big)dyds-
$$
$$
-\int\limits_a^{t}
\int\limits_{\Omega}W(t,s,x,y)f\Big(\big(S_{\tau}^\varphi
u_2\big)(s,x,y)\Big)dyds\Big|\leq
$$
$$
\leq\sup\limits_{t\in[a+\gamma,a+\gamma+\delta],
x\in\Omega}\Big|\int\limits_{a+\gamma}^{t}
\int\limits_{\Omega}W(t,s,x,y)\Bigg(f\Big(\big(S_{\tau}^\varphi
u_1\big)(s,x,y)\Big)-\ \ \ \ \
 \ \ \ \ \ \  \  \ \ \ \ \   \ $$ $$\ \ \ \ \  \ \  \ \ \ \ \   \ \ \ \ \  \ \ \ \ \ \ \ \  \
  \ \ \ \ \
  \  \ \ \ \ \  \ \ \ \ \ \ \ \ \ -f\Big(\big(S_{\tau}^\varphi
u_2\big)(s,x,y)\Big)\Bigg)dyds\Big|\leq
$$
$$
\leq\sup\limits_{t\in[a+\gamma,a+\gamma+\delta],
x\in\Omega}\int\limits_{a+\gamma}^{t}
\int\limits_{\Omega}\Big|W(t,s,x,y)\Big|
\widetilde{f}_r\|u_1-u_2\|_{C([a,b],BC(\Omega, R^n))}dyds\leq
$$
$$
\leq\sup\limits_{t\in[a+\gamma,a+\gamma+\delta],
x\in\Omega}\int\limits_{a+\gamma}^{t}
\int\limits_{\Omega}\Big|W(t,s,x,y)\Big|\widetilde{f}_r
dyds\|u_1-u_2\|_{C([a,b],BC(\Omega, R^n))}\leq$$$$\leq q
\|u_1-u_2\|_{C([a,b],BC(\Omega, R^n))}.
$$
Here
$$q=\widetilde{f}_r\sup\limits_{t\in[a+\gamma,a+\gamma+\delta],
x\in\Omega}\int\limits_{a+\gamma}^{t}
\int\limits_{\Omega}\Big|W(t,s,x,y)\Big| dyds.$$ Thus, we can always
find $\delta>0$ such that $q\leq q_0$. Hence, the property
$\mathfrak{q}_2)$ for the mapping $F$, given by $(3.2)$, holds true.
The verification of the property $\mathfrak{q}_1)$ is analogous.
Taking into account  Remark 2.2, we prove the theorem. \end{proof}

\begin{remark}If in the Theorem 3.1 condition
$\widetilde{f}_r=\widetilde{f}$ is independent of $r$ (as e.g. in
classical neural field models, where $0\leq f(u) \leq 1$), then
according to   Remark 2.1 we will get a global solution to the Eq.
$(3.1)$. In this case, if we  take $\tau(t,x,y)\equiv0$,
 Theorem 3.1 becomes
analogous to the results concerning solvability of the Amari model
obtained by Potthast \textit{et al}. \cite{PG}.
\end{remark}

\null

\begin{remark}If in  Theorem 3.1 the condition
$\widetilde{f}_r=\widetilde{f}$ is independent of $r$, Theorem 3.1
can be compared to the theorem on solvability
 of  Eq. $(1.3)$ in $C([a,b],L^2(\Omega, R^n))$ for any $b>a$ proved in Faye \textit{et al}. \cite{Faye} Here
 we obtained the same result for the more general model $(3.1)$ in $C([a,b],BC(\Omega, R^n))$.
 We note that in case when the delay $\tau(t,x,y)=\tau(x,y)$ is independent of $t$, it is possible
 to prove Theorem 3.1 for the space $C([a,b],L^2(\Omega, R^n))$ using
 our
 technique as well thus getting the main theoretical result of
 \cite{Faye}.
\end{remark}

\null

Note that the remarks 3 and 4  on maximally extended solutions are
valid for the problem $(3.1)$ as well.

It is also worth mentioning that our approach to delayed
functional-differential equations is based on the idea to include
the prehistory function in the inner superposition operator. It
allows us to consider the operator equation $(2.1)$ with the
  operator $(3.2)$ defined on $[a,b]$ instead of $(-\infty,b]$. The
same approach to   functional-differential equations with delay was
implemented e.g. in  \cite{Ber2005},  \cite{Ber2011}.

\null

Next we complete the study of wellposedness of the problem $(3.1)$
by
 investigating continuous dependence of solutions to the associated
problem

 \begin{eqnarray}
\begin{array}{c}\displaystyle
u(t,x)=\varphi_{\lambda}(a,x)+\int\limits_a^t
 \int\limits_{\Omega}W_{\lambda}(t,s,x,y)
f_{\lambda}(u(s-\tau_{\lambda}(s,x,y),y))dyds,\\
 t\in[a,\infty), x\in\Omega;   \\
u(\xi,x)= \varphi_{\lambda}(\xi,x),\ \xi\leq a,
x\in\Omega\end{array}
 \end{eqnarray}
on a parameter  $\lambda\in \Lambda$.

The assumptions $(A_\lambda1)$ -- $(A_\lambda6)$ imposed on the
functions in the model $(3.3)$ for each $\lambda\in \Lambda$ repeat
the  assumptions  $(A1)$ -- $(A6)$, respectively.

\null

We will naturally apply Definition 3.1 to the model $(3.3)$ at each
$\lambda\in \Lambda$.

 The following theorem gives conditions that guarantee wellposedness of the
 problem $(3.3)$.

 \null

\textbf{  Theorem 3.2.} \textit{Let the assumptions $(A_\lambda1)$
-- $(A_\lambda6)$ hold true. Assume that the following  conditions
are satisfied:}

\vspace{10pt}

 1) \textit{There is a neighborhood $U_0$ of
$\lambda_0$ such that for any   $r>0$ there exists
$\widetilde{f}_r\in R$ (independent of $\lambda\in U_0$) such that
for which $|f_\lambda(u_1)-f_\lambda(u_2)|\leq
\widetilde{f}_r|u_1-u_2|$ for all $u_1, u_2 \in R^n$, $|u_1|\leq r$,
$|u_2|\leq r$.}

\vspace{10pt}

\textit{For any $\{\lambda_i\}\subset\Lambda$,
$\lambda_i\to\lambda_0$ it holds true that:}

\vspace{10pt}

2) \textit{For any $b>a$,
 $$\sup\limits_{t\in[a,b],\ x\in\Omega}\Big|\int\limits_{a}^t\int\limits_{\Omega}W_{\lambda_i}(t,s,x,y)dyds-
\int\limits_{a}^t\int\limits_{\Omega}W_{\lambda_0}(t,s,x,y)dyds\Big|\to0;$$}

\vspace{10pt}

 3) \textit{For any $b>a$,  if  $|u_i(\cdot,\cdot)-u(\cdot,\cdot)|\to 0$ in measure on $[a,b]\times\Omega$
  as $i\to\infty$, then
$|f_{\lambda_i}(u_i(\cdot,\cdot))-f_{\lambda_0}(u(\cdot,\cdot))|\to
0$ in measure on $[a,b]\times\Omega$
  as $i\to\infty$;}

\vspace{10pt}

 4) \textit{For any $b>a$,
$\sup\limits_{x\in\Omega}|\tau_{\lambda_i}(\cdot,x,\cdot)-\tau_{\lambda_0}(\cdot,x,\cdot
)|\to 0$ in measure on $[a, b]\times\Omega$;}

\vspace{10pt}

5) $\|\varphi_{\lambda_i} - \varphi_{\lambda_0}\|_{C((-\infty,a],
BC(\Omega,R^n))}\to0$.

\vspace{10pt}

  \textit{ Then there is a  neighborhood $U$ of $\lambda_0$,
 such that for each element $\lambda\in U$, Eq.
$(3.3)$ has a unique global or maximally extended solution, and each
local solution is a restriction of this  solution.
  Moreover, if
  at $\lambda=\lambda_0$
  Eq. $(3.3)$
 has a local solution $u_{0\gamma}$ defined on
 $[a,a{+}\gamma]\times\Omega$, then
  for any $\{\lambda_i\}\subset\Lambda$, $\lambda_i\to\lambda_0$ one
  can find number $I$ such that for all $i>I$ Eq. $(3.3)$  has a local
solution $u_{\gamma}=u_{\gamma}(\lambda_i)$ defined on $[a,
a{+}\gamma]\times\Omega$ and $\|u_{\gamma}(\lambda_i)-
u_{0\gamma}\|_{C([a,a+\gamma],BC (\Omega,R^n))}\rightarrow0$.}

\begin{proof} Choose an
 arbitrary $b\in(a,\infty)$. In order to use  Theorem 2.1, we need to bring the Eq. $(3.3)$ to the
form   $u(t,\cdot)=(F(u,\lambda))(t), \ t\in[a,b]$. Using the same
technique as in the proof of Theorem 3.1 and the corresponding
assumptions $(A_\lambda1)$ - $(A_\lambda6)$, we get here
$$F(\cdot,\lambda):C([a,b],BC(\Omega,R^n))\to
C([a,b],BC(\Omega,R^n)),
$$
$$
(F(u,\lambda))(t,x)=\varphi_\lambda(a,x)+\int\limits_a^t
 \int\limits_{\Omega}W_\lambda(t,s,x,y)f_\lambda\Big(\big(
S_{\tau_\lambda}^{\varphi_\lambda} u\big)(s,x,y)\Big)dyds, \
$$
$$t\in[a,b],\ x\in \Omega,$$
$$
\big(S_{\tau_\lambda}^{\varphi_\lambda} u\big)(t,x,y)=\left\{
 \begin{array}{ccl}
  \varphi_\lambda(t-{\tau_\lambda}(t,x,y),y), \!\!\!&\mbox{if}\!\!\!&~t-{\tau_\lambda}(t,x,y)<a;\\
  u(t-{\tau_\lambda}(t,x,y),y),\!\!\!&\mbox{if}\!\!\!&~t-{\tau_\lambda}(t,x,y)\geq a\\
 \end{array}
\right.
$$
for all $\lambda \in \Lambda$.

 The condition 1) of   this theorem
allows us to verify the assumption 1) of  Theorem 2.1   for each
$\lambda \in U_0$ by the same procedure as we used in the proof of
Theorem~2. So, we only need to verify the condition 2) of Theorem~1.

Choose an arbitrary $u\in C([a,b], BC(\Omega,R^n))$. Let
$\|u_i-u\|_{Y}\to0$, i.e, $\|u_i-u\|_{C([a,b],
BC(\Omega,R^n))}\to0$, $i\to\infty$, and $\lambda\to\lambda_0$.

We have the following estimates:
$$
|(S_{\tau_\lambda}^{\varphi_\lambda}
u_i)(t,x,y)-(S_{\tau_{\lambda_0}}^{\varphi_{\lambda_0}}
u)(t,x,y)|\leq|(S_{\tau_\lambda}^{\varphi_\lambda}
u_i)(t,x,y)-(S_{\tau_\lambda}^{\varphi_\lambda} u)(t,x,y)|+
$$
$$
+|(S_{\tau_\lambda}^{\varphi_\lambda}
u)(t,x,y)-(S_{\tau_{\lambda_0}}^{\varphi_\lambda}
u)(t,x,y)|+|(S_{\tau_{\lambda_0}}^{\varphi_\lambda}
u)(t,x,y)-(S_{\tau_{\lambda_0}}^{\varphi_{\lambda_0}} u)(t,x,y)|.
$$

If $\lambda\to\lambda_0$, then the first term on the right-hand side
of this inequality goes to $0$ uniformly  as $\|u_i-u\|_{C([a,b],
BC(\Omega,R^n))}\to0$.
   By the virtue of the condition 4), the second term on
the right-hand side goes to $0$ in measure on $([a,b]
\times\Omega)$, uniformly in $x\in\Omega$, as $\lambda\to\lambda_0$.
The third term on the right-hand side of the inequality goes to $0$
uniformly when $\lambda\to\lambda_0$ as the condition 5) holds true.
Thus,  we have
$$|(S_{\tau_\lambda}^{\varphi_\lambda}
u_i)(\cdot,x,\cdot)-(S_{\tau_{\lambda_0}}^{\varphi_{\lambda_0}}
u)(\cdot,x,\cdot)|\to 0$$ in measure, uniformly in $x\in\Omega$, as
$\|u_i-u\|_{C([a,b], BC(\Omega,R^n))}\to0$ and
$\lambda\to\lambda_0$.

  Using this convergence, we can make the following estimates
$$
\sup\limits_{t\in[a,b],\ x\in\Omega}\Big|\int\limits_a^t
 \int\limits_{\Omega}W_{\lambda}(t,s,x,y)f_\lambda\Big(\big(S^{\varphi_{\lambda}}_{\tau_{\lambda}}
u_i\big)(s,x,y)\Big)dyds- \ \ \ \ \ \ \ \ \ \ \ \ \ \ \ \ \ \ \ \ \
\
$$$$\ \ \ \ \ \ \ \ \ \ \ \ \
\ \ \ \ \ \ \ \-\int\limits_a^t
 \int\limits_{\Omega}W_{\lambda_0}(t,s,x,y)f_{\lambda_0}\Big(\big(S^{\varphi_{\lambda_0}}_{\tau_{\lambda_0}}
u\big)(s,x,y)\Big)dyds\Big|\leq
$$
$$
\sup\limits_{t\in[a,b],\ x\in\Omega}\Big|\int\limits_a^t
 \int\limits_{\Omega}W_{\lambda}(t,s,x,y)f_\lambda\Big(\big(S^{\varphi_{\lambda}}_{\tau_{\lambda}}
u_i\big)(s,x,y)\Big)dyds- \ \ \ \ \ \ \ \ \ \ \ \ \ \ \ \ \ \ \ \ \
\
$$$$\ \ \ \ \ \ \ \ \ \ \ \ \
\ \ \ \ \ \ \ \-\int\limits_a^t
 \int\limits_{\Omega}W_{\lambda}(t,s,x,y)f_{\lambda_0}\Big(\big(S^{\varphi_{\lambda_0}}_{\tau_{\lambda_0}}
u\big)(s,x,y)\Big)dyds\Big|+
$$
$$
+\sup\limits_{t\in[a,b],\ x\in\Omega}\Big|\int\limits_a^t
 \int\limits_{\Omega}W_{\lambda}(t,s,x,y)f_{\lambda_0}\Big(\big(S^{\varphi_{\lambda_0}}_{\tau_{\lambda_0}}
u\big)(s,x,y)\Big)dyds- \ \ \ \ \ \ \ \ \ \ \ \ \ \ \ \ \ \ \ \ \ \
$$$$\ \ \ \ \ \ \ \ \ \ \ \ \
\ \ \ \ \ \ \ \-\int\limits_a^t
 \int\limits_{\Omega}W_{\lambda_0}(t,s,x,y)f_{\lambda_0}\Big(\big(S^{\varphi_{\lambda_0}}_{\tau_{\lambda_0}}
u\big)(s,x,y)\Big)dyds\Big|.
$$
Taking into account the condition 3), we conclude that the first
term on the right-hand side of the inequality goes to $0$ as
$\lambda\to\lambda_0$. The second   term  on the right-hand side of
the inequality goes to $0$ by the virtue of the condition     2) as
$\lambda\to\lambda_0$.

Thus,  the condition 2) of   Theorem 2.1 is satisfied and Theorem
3.2 is proved. \end{proof}

\null

We emphasize here that our aim was to formulate the assumptions on
the functions involved in the model $(3.3)$ (see conditions 2) -- 5)
of Theorem 3.2) as general as it possible. Of course, we can
strengthen these assumptions in order to make them more conventional
e.g. in the following way.

\null

\begin{remark}If the estimate in the assumption $(A_\lambda3)$ holds true
uniformly with respect to $\lambda\in\Lambda $, then it is possible
to get the conclusion of Theorem 3.2 by claiming that for any $b>a$
the functions
 $$
 W_{(\cdot)}:\Lambda\times  [a,b] \times [a,b] \times \Omega \times
 \Omega \to R^n,
 $$
 $$
 f_{(\cdot)}: \Lambda\times R^n\to R^n,
 $$
$$
\tau_{(\cdot)}:\Lambda  \times [a,b] \times \Omega \times  \Omega
\to [0,\infty),
$$
$$
\varphi_{(\cdot)}: \Lambda \times(-\infty,b]\times \Omega\to R^n$$
 are continuous instead of claiming the conditions $2)$ --
$5)$ of Theorem 3.2.
\end{remark}

\null

We  now  consider two important special cases of the model $(3.3)$.

As the neural field theory studies processes in cortical tissue, it
is realistic to assume that $\Omega$ is bounded (see e.g.
 \cite{Faye}). The following remark represents the result, analogous
to Theorem 3.2 for this case.

\null

\begin{remark}If $\Omega$ is bounded,   we can  substitute $(A_\lambda6)$
by

\vspace{10pt}

 $(A_\lambda^\ast6)$ For any $a^\ast<a$ and each
$\varphi_\lambda\in C([a^\ast,a],C_{0}(\Omega,R^n))$, $\lambda\in\Lambda$. \\
In order to get the conclusion of   Theorem 3.2, we need the
following conditions instead of the conditions $3)$, $4)$, and $5)$,
respectively:

\vspace{10pt}

 For any $\{\lambda_i\}\subset\Lambda$,
$\lambda_i\to\lambda_0$ it holds true that:

\vspace{10pt}

 3$^\ast$)  For any  $u\in R^n$ we have
$|f_{\lambda_i}(u)-f_{\lambda_0}(u)|\to 0$;   \vspace{10pt}

 4$^\ast$)  For all $x\in\Omega$,
$|\tau_{\lambda_i}(\cdot,x,\cdot)-\tau_{\lambda_0}(\cdot,x,\cdot)|\to
0$ in measure on $[a, b]\times\Omega$;

\vspace{10pt}

5$^\ast$) For any  $a^\ast<a$ and all
$(t,x)\in[a^\ast,a]\times\Omega$, $|\varphi_{\lambda_i}(t,x) -
\varphi_{\lambda_0}(t,x)|\to0$.

\end{remark}

\null

Proof of the statement in Remark 3.4 is given in Appendix A.

\null

In neural field modeling special attention is paid to spatially
localized solutions, so-called "bumps". If  $\Omega$ is unbounded,
but the  solution to $(3.3)$ is spatially localized, we can relax
  Theorem 3.2 conditions in the following way.

\null

\begin{remark}If we replace $(A_\lambda6)$ by

 \vspace{10pt}

  $(A_\lambda^\prime6)$  For each $\lambda\in \Lambda$, the
prehistory function $\varphi_\lambda\in
C((-\infty,a],C_{0}(\Omega,R^n))$;\\
 and impose the additional condition, corresponding to localization
 in the spatial variable,

 \vspace{10pt}

$(A_\lambda^\prime7)$  For each $\lambda\in \Lambda$ and any $b>a$,
$\lim\limits_{|x|\to\infty}\big| W_\lambda(t,s,x,y)\big|= 0$ for all
 $(t,s,y)\in [a,b]\times[a,b]\times\Omega$,
then, in order to get the conclusion of   Theorem 3.2 holds true for
spatially localized solutions,  we need the following conditions,
instead of  2), 3), 4), and 5)respectively:

\vspace{10pt}

 For any $\{\lambda_i\}\subset\Lambda$,
$\lambda_i\to\lambda_0$ it holds true that:

\vspace{10pt}

2$^\prime$)  For any $b>a$, $r>0$, and each $t\in[a,b]$, $x\in
\Omega, \ |x|\leq r$ it holds true that
$$\Big|\int\limits_{a}^t\int\limits_{\Omega}\Big(W_{\lambda_i}(t,s,x,y)dy-
W_{\lambda_0}(t,s,x,y)\Big)dyds\Big|\to 0;$$

\vspace{10pt}

 3$^\prime$)  For any  $u\in R^n$ we have
$|f_{\lambda_i}(u)-f_{\lambda_0}(u)|\to 0$;

\vspace{10pt}

 4$^\prime$)  For all $x\in\Omega$,
$|\tau_{\lambda_i}(\cdot,x,\cdot)-\tau_{\lambda_0}(\cdot,x,\cdot)|\to
0$ in measure on $[a, b]\times\Omega$;

\vspace{10pt}

5$^\prime$)  For any $(t,x)\in(-\infty,a]\times\Omega$,
$|\varphi_{\lambda_i}(t,x) - \varphi_{\lambda_0}(t,x)|\to0$.
\end{remark}

\null

 Proof of the statement in Remark 3.5 is given in Appendix B.

\null

As Theorems 2 and 3 are valid for each $a\in R$ in the model
$(3.1)$, it is natural to address the question, what happens in the
case when $a=-\infty$ (i.e., when $(3.1)$ becomes $(1.7)$).

\null

\begin{remark}Solution to  $(1.7)$ is not
necessarily unique.
\end{remark}

\null

The following example illustrates this fact.

\null

 \textbf{Example 3.1.}
Consider the equation
$$u(t,x)=\int\limits_{-\infty}^t
\int\limits_{R}\exp(-s){\omega}(x) u(s,y)dyds, \ t\in R, \ x\in R $$
with some Gaussian function ${\omega}$.
 Define the function $u\in C((-\infty,b],BC(R,R))$ as follows: $$u(t,x)=v(t)\omega(x), $$ where
 $$v(t)=V \exp \big( -\exp(-t) \big) ,  V\in R,$$ is a solution to $$\dot{v}(t)=\exp(-t)v(t),$$ satisfying the
 property $v(t)\to0$ as $t\to-\infty$.
Thus, for any $V\in R$ we get
  a solution to
 $(1.7)$ which belongs to $ C((-\infty,b],BC(R,R))$.

 \null

Nevertheless, it is possible to find conditions, which guarantee
wellposedness of the model $(1.7)$. The last part of the present
paper is devoted to this problem. We have the following assumptions
on the functions involved:

\vspace{10pt}

$(\mathcal{A}1)$  For any $a,b\in R$, $a<b$, $t \in [a, b]$,
$x\in\Omega$, the function $W(t,\cdot,x,\cdot):[a,b]\times\Omega\to
R^n$ is measurable.

\vspace{5pt}

$(\mathcal{A}2)$ For any  $a,b\in R$, $a<b$, at almost all $(s,y)
\in [a, b]\times\Omega$, the function
$W(\cdot,s,\cdot,y):(-\infty,b]\times \Omega\to R^n$ is uniformly
continuous.

 \vspace{5pt}

$(\mathcal{A}3)$ For any  $b\in R$, $t \in (-\infty,b]$,
$\int\limits_{\Omega}\sup\limits_{x\in\Omega}\big| W(t,s,x,y)\Big|
dy = G(s),$
 where $G\in  L^1((-\infty,b],\mu,R^n)$.

\vspace{5pt}

Assumptions $(\mathcal{A}4)$ and
  $(\mathcal{A}5)$ are the same as the corresponding assumptions $(A4)$ and
  $(A5)$.

 \vspace{10pt}

Now, we need to give the definitions of local, maximally extended
and global solutions  to  Eq. $(1.7)$.

\null

 \textbf{Definition 3.2.}
     We define \textit{a local solution} to Eq. $(1.7)$  on $(-\infty,\gamma]\times\Omega$,
 $\gamma\in R$,
to be a function $u_\gamma\in  C((-\infty, \gamma],BC(\Omega,R^n))$
that satisfies the equation $(1.7)$ on $(-\infty,
\gamma]\times\Omega$. We define a \textit{maximally extended
solution} to Eq. $(1.7)$  on $(-\infty, {\zeta})\times\Omega$,
${\zeta}\in R$ to be a function $u_{\zeta}:(-\infty,
{\zeta})\times\Omega\rightarrow
 {R}^n$, whose restriction $u_\gamma$ to $(-\infty,  \gamma]\times\Omega$
 is a local solution to Eq. $(1.7)$ for any $\gamma<{\zeta}$ and
$\lim\limits_{\gamma\rightarrow{\zeta}-0}\|u_\gamma\|_{C((-\infty,
\gamma],BC(\Omega,R^n))}=\infty$. We define a \textit{global
solution} to Eq. $(1.7)$ to be a function $u:R\times\Omega\to R^n$,
whose restriction $u_\gamma$ to $(-\infty,  \gamma]\times\Omega$ is
its local solution for any $\gamma\in R$.

\null

 \textbf{  Theorem 3.3.}
\textit{Let the assumptions  $(\mathcal{A}4)$ --
  $(\mathcal{A}5)$ hold true. If for any $r>0$ there exists
$\widetilde{f}_r\in R$ such that for all $u_1, u_2 \in R^n$,
$|u_1|\leq r$, $|u_2|\leq r$, we have $|f(u_1)-f(u_2)|\leq
\widetilde{f}_r|u_1-u_2|$, then
 Eq.
$(1.7)$ has a unique global or maximally extended solution and each
local solution is a restriction of this global or
 maximally extended solution (all types of solutions are meant in the sense of Definition 3.2).}

 \begin{proof}
 First,  we  prove existence of a unique local
 solution to $(1.7)$.
Choose arbitrary $b\in R$. Using the same estimation technique as in
the proof ot Theorem 3.1 and the corresponding assumptions
$(\mathcal{A}1)$ --
  $(\mathcal{A}5)$, we rewrite   Eq. $(1.7)$
as the operator equation $u(t,\cdot)=(Fu)(t)$, and consider it on
$(-\infty,b]$, where
$$F: C((-\infty,b],BC(\Omega,R^n))\to
 C((-\infty,b],BC(\Omega,R^n)),
$$
$$
(Fu)(t,x)=\int\limits_{-\infty}^t
 \int\limits_{\Omega}W(t,s,x,y)f\Big(u(s-{\tau}(s,x,y),y)\Big)dyds,
\ t\in[a,b],\ x\in \Omega.
$$

Choose arbitrary $q_0<1$, $r>0$,
$\|u_1\|_{C((-\infty,b],BC(\Omega,R^n))}\leq r$ and
$\|u_2\|_{C((-\infty,b],BC(\Omega,R^n))}\leq r$. In order to prove
existence of a unique local solution to $(1.7)$ using the Banach
fixed point theorem, we need to find $\delta\in R$ such that
$$\max\limits_{t\in
(-\infty,\delta]}\|(Fu_1)(t)-(Fu_2)(t)\|_{BC(\Omega,R^n)}\leq q_0
\max\limits_{t\in
(-\infty,\delta]}\|(u_1)(t)-(u_2)(t)\|_{BC(\Omega,R^n)}.$$

For any $\delta<b$, we get the estimates
$$
\sup\limits_{t\in(-\infty,\delta],
x\in\Omega}\Big|\int\limits_{-\infty}^{t}
\int\limits_{\Omega}W(t,s,x,y)f\Big(u_1(s-{\tau}(s,x,y),y)\Big)dyds-
$$
$$
-\int\limits_{-\infty}^{t}
 \int\limits_{\Omega}W(t,s,x,y)f\Big(u_2(s-{\tau}(s,x,y),y)\Big)dyds\Big|\leq
$$
$$
\leq\sup\limits_{t\in(-\infty,\delta],
x\in\Omega}\Big|\int\limits_{-\infty}^{t}
 \int\limits_{\Omega}W(t,s,x,y)\Bigg(f\Big(u_1(s-{\tau}(s,x,y),y)\Big)-\
\ \ \ \
 \ \ \ \ \ \  \  \ \ \ \ \   \ $$ $$\ \ \ \ \  \ \  \ \ \ \ \   \ \ \ \ \  \ \ \ \ \ \ \ \  \
  \ \ \ \ \
  \  \ \ \ \ \  \ \ \ \ \ \ \ \ \ -f\Big(u_2(s-{\tau}(s,x,y),y)\Bigg)dyds\Big|\leq
$$
$$
\leq\sup\limits_{t\in(-\infty,\delta],
x\in\Omega}\int\limits_{-\infty}^{t}
 \int\limits_{\Omega}\Big|W(t,s,x,y)\Big|\widetilde{f}_r
dyds\|u_1-u_2\|_{BC((-\infty,\delta]\times \Omega,R^n)}\leq$$$$\leq
q \|u_1-u_2\|_{BC((-\infty,\delta]\times \Omega,R^n)}.
$$
Here
$$q=\widetilde{f}_r\sup\limits_{t\in(-\infty,\delta], x\in\Omega}\int\limits_{-\infty}^{t}
 \int\limits_{\Omega}\Big|W(t,s,x,y)\Big| dyds.$$ Using the
assumption $(\mathcal{A}3)$, we can  find $\delta>0$ such that
$q\leq q_0$. Thus, the equation $(1.7)$ has a unique local solution,
defined on $(-\infty,\delta]\times\Omega$. Now, regarding this
solution as a prehistory
 function for the model $(3.1)$ and taking $a=\delta$, we use Theorem 3.1 and obtain the
 conclusion of the theorem.
\end{proof}

In order to approach the problem of wellposedness of $(1.7)$, we
consider its parameterized version:

$$
u(t,x)= \int\limits_{-\infty}^t
\int\limits_{\Omega}W_{\lambda}(t,s,x,y)
f_{\lambda}(u(s-\tau_{\lambda}(s,x,y),y))dyds,\  $$$$
 t\in R, x\in \Omega, \eqno(3.4)$$
with a parameter  $\lambda\in \Lambda$.

For each $\lambda\in \Lambda$, the assumptions
$(\mathcal{A}_\lambda1)$ -- $(\mathcal{A}_\lambda5)$, imposed on the
functions involved in the model $(3.4)$ repeat the  assumptions
$(\mathcal{A}1)$ -- $(\mathcal{A}5)$, respectively.

\null

 At each $\lambda\in\Lambda$ we define the types of
solutions to $(3.4)$ according to Definition~3.2.

\null

\textbf{  Theorem 3.4.} \textit{Let the assumptions
$(\mathcal{A}_\lambda1)$ --
  $(\mathcal{A}_\lambda5)$ hold true. Assume that the following  conditions are
satisfied:}

\vspace{10pt}

 1) \textit{There is a neighborhood $U_0$ of
$\lambda_0$ such that for any for any $r>0$ there exists
$\widetilde{f}_r\in R$ (independent of $\lambda\in U_0$),  for which
$|f_\lambda(u_1)-f_\lambda(u_2)|\leq \widetilde{f}_r|u_1-u_2|$  for
all $u_1, u_2 \in R^n$, $|u_1|\leq r$, $|u_2|\leq r$;}

\vspace{10pt}

\textit{For any $\{\lambda_i\}\subset\Lambda$,
$\lambda_i\to\lambda_0$ it holds true that:}

\vspace{10pt}

2) \textit{ For any $b\in R$,
 $\sup\limits_{(-\infty,b],\ x\in\Omega}\Big|\int\limits_{-\infty}^t\int\limits_{\Omega}W_{\lambda_i}(t,s,x,y)dy-
\int\limits_{-\infty}^t\int\limits_{\Omega}W_{\lambda_0}(t,s,x,y)dy\Big|\to
0;$}

\vspace{10pt}

 3) \textit{For any $b\in R$,  if  $|u_i(\cdot,\cdot)-u(\cdot,\cdot)|\to 0$ in measure on $(-\infty,b]\times\Omega$
  as $i\to\infty$, then
$|f_{\lambda_i}(u_i(\cdot,\cdot))-f_{\lambda_0}(u(\cdot,\cdot))|\to0$
in measure on $(-\infty,b]\times\Omega$
  as $i\to\infty$;}

\vspace{10pt}

4) \textit{For any $b\in R$,
$\sup\limits_{x\in\Omega}|\tau_{\lambda_i}(\cdot,x,\cdot)-\tau_{\lambda_0}(\cdot,x,\cdot)|\to
0$ in measure on $(-\infty,b]\times\Omega$; }

\vspace{10pt}

  \textit{ Then there is a  neighborhood $U$ of $\lambda_0$,
 such that for each  $\lambda\in U$, Eq.
$(3.4)$ has a unique global or maximally extended solution, and each
local solution is a restriction of this  solution.
  Moreover, if
  at $\lambda=\lambda_0$
  Eq. $(3.4)$
 has a local solution $u_{0\gamma}$ defined on
 $(-\infty,\gamma]\times\Omega$, then
  for any $\{\lambda_i\}\subset\Lambda$, $\lambda_i\to\lambda_0$ one
  can find number $I$ such that for all $i>I$ Eq. $(3.4)$  has a local
solution $u_{\gamma}=u_{\gamma}(\lambda_i)$ defined on
$(-\infty,\gamma]\times\Omega$ and $\|u_{\gamma}(\lambda)-
u_{0\gamma}\|_{C((-\infty,\gamma],BC (\Omega,R^n))}\rightarrow0$ as
$\lambda\to\lambda_0$.}

\begin{proof} Choose an
 arbitrary $b\in R$. Consider the following operator equation
   $$u(t,\cdot)=(F(u,\lambda))(t), \ t\in(-\infty,b],$$ where at each
$\lambda \in \Lambda$, by the virtue of the assumptions
$(\mathcal{A}_\lambda1)$ -- $(\mathcal{A}_\lambda5)$,
$$
F(\cdot,\lambda):C((-\infty,b],BC(\Omega,R^n))\to
C((-\infty,b],BC(\Omega,R^n)),
$$
$$
(F(u,\lambda))(t,x)= \int\limits_{-\infty}^t
\int\limits_{\Omega}W_\lambda(t,s,x,y)f_\lambda
\Big(u(t-{\tau_\lambda}(t,x,y),y)\Big)dyds, .
$$
$$ t\in(-\infty,b],\
x\in\Omega$$

  Note that by Theorem 3.3  we have a
unique solution to  Eq. $(3.4)$ defined on
$(-\infty,\delta]\times\Omega$ for each $\lambda \in U_0$. We need
to prove continuous dependence of these solutions on $\lambda$.
First, we prove that the operator $F$ is continuous in
$(u,\lambda_0)$ for any fixed $u\in C((-\infty,b],BC (\Omega,R^n))$.

Choose an arbitrary $u\in C((-\infty,b],BC (\Omega,R^n))$. Let
  $\|u_i-u\|_{C((-\infty,b],BC (\Omega,R^n))}\to0$, $i\to\infty$, and $\lambda\to\lambda_0$.

We have the following estimates:
$$
|u_i(t-{\tau_\lambda}(t,x,y),y)-u(t-{\tau_{\lambda_0}}(t,x,y),y)|\leq
$$
$$\leq|u_i(t-{\tau_\lambda}(t,x,y),y)-u(t-{\tau_{\lambda}}(t,x,y),y)|+$$
$$
+|u(t-{\tau_{\lambda}}(t,x,y),y)-u(t-{\tau_{\lambda_0}}(t,x,y),y)|.
$$

If $\lambda\to\lambda_0$, then the first term on the right-hand side
of this inequality goes to $0$ uniformly  as $\|u_i-u\|_{
C((-\infty,b],BC (\Omega,R^n))}\to0$.
   By  virtue of the condition 4), the second term on
the right-hand side goes to $0$ in measure on $((-\infty,b]
\times\Omega)$, uniformly in $x\in\Omega$, as $\lambda\to\lambda_0$.
So,
$$|u_i(\cdot-{\tau_\lambda}(\cdot,x,\cdot),\cdot)-u(t-{\tau_{\lambda_0}}(\cdot,x,\cdot),\cdot)| \to 0$$ in measure, uniformly in $x\in\Omega$, as
$\|u_i-u\|_{C((-\infty,b],BC (\Omega,R^n))}\to0$ and
$\lambda\to\lambda_0$.

  Using this convergence, we obtain
$$
\sup\limits_{t\in(-\infty,b],\
x\in\Omega}\Big|\int\limits_{-\infty}^t
 \int\limits_{\Omega}W_{\lambda}(t,s,x,y)f_\lambda
\Big(u_i(t-{\tau_\lambda}(s,x,y),y)\Big)dyds- \ \ \ \ \ \ \ \ \ \ \
\ \ \ \ \ \ \ \ \ \ \
$$$$\ \ \ \ \ \ \ \ \ \ \ \ \
\ \ \ \ \ \ \ \-\int\limits_{-\infty}^t
 \int\limits_{\Omega}W_{\lambda_0}(t,s,x,y)f_{\lambda_0}\Big
(u(t-{\tau_{\lambda_0}}(s,x,y),y)\Big)dyds\Big|\leq
$$
$$
\sup\limits_{t\in(-\infty,b],\
x\in\Omega}\Big|\int\limits_{-\infty}^t
 \int\limits_{\Omega}W_{\lambda}(t,s,x,y)f_\lambda\Big(u_i(t-{\tau_\lambda}(s,x,y),y)\Big)dyds-
\ \ \ \ \ \ \ \ \ \ \ \ \ \ \ \ \ \ \ \ \ \
$$$$\ \ \ \ \ \ \ \ \ \ \ \ \
\ \ \ \ \ \ \ \-\int\limits_{-\infty}^t
 \int\limits_{\Omega}W_{\lambda}(t,s,x,y)f_{\lambda_0}
\Big(u(t-{\tau_{\lambda_0}}(s,x,y),y)\Big)dyds\Big|+
$$
$$
+\sup\limits_{t\in(-\infty,b],\
x\in\Omega}\Big|\int\limits_{-\infty}^t
 \int\limits_{\Omega}W_{\lambda}(t,s,x,y)f_{\lambda_0}
\Big(u(t-{\tau_{\lambda_0}}(s,x,y),y)\Big)dyds- \ \ \ \ \ \ \ \ \ \
\ \ \ \ \ \ \ \ \ \ \ \
$$$$\ \ \ \ \ \ \ \ \ \ \ \ \
\ \ \ \ \ \ \ \-\int\limits_{-\infty}^t
 \int\limits_{\Omega}W_{\lambda_0}(t,s,x,y)f_{\lambda_0}
\Big(u(t-{\tau_{\lambda_0}}(s,x,y),y)\Big)dyds\Big|.
$$
Taking into account the condition 3), we conclude that the first
term on the right-hand side of the inequality goes to $0$ as
$\lambda\to\lambda_0$. The second  term  on the right-hand side of
the inequality goes to $0$ by the virtue of the condition    2) as
$\lambda\to\lambda_0$. Thus, the operator $F$ is continuous in
$(u,\lambda_0)$ for any chosen $u\in C((-\infty,b],BC
(\Omega,R^n))$. Using this fact, for any $\varepsilon>0$ we can find
such  $\varepsilon_1>0$ and neighborhood $U_1$ of $\lambda_0$, that
$$\|F(\mathfrak{u}_{\delta},\lambda)-F(u_{0\delta},\lambda)\|_{C((-\infty,b],BC (\Omega,R^n))}\leq \varepsilon$$ for all $\lambda\in U_1$ and any
$\mathfrak{u}_{\delta}\in C((-\infty,\delta],BC(\Omega,R^n))$,
satisfying the estimate
$$\|\mathfrak{u}_{\delta}-u_{0\delta}\|_{C((-\infty,\delta],BC(\Omega,R^n))}\leq
\varepsilon_1.$$

As the mapping $F(\cdot,\lambda)$ is contracting with the constant
$q_0<1$ (see Theorem 3.3) for any  $\lambda\in U_0$, for any $m=1,2,
\ldots$ we have
$$\|F^m(u_{0\delta},\lambda)-u_{0\delta}\|_{C((-\infty,\delta],BC(\Omega,R^n))}\leq $$$$ \leq\|F^m(u_{0\delta},\lambda)-
F^{m-1}(u_{0\delta},\lambda)\|_{C((-\infty,\delta],BC(\Omega,R^n))}+
\ldots$$
$$\ldots+\|F(u_{0\delta},\lambda)-u_{0\delta}\|_{C((-\infty,\delta],BC(\Omega,R^n))}\leq $$$$ \leq(q_0^{m-1}+\ldots+q_0+1)
(1-q_0)\varepsilon\leq\varepsilon.$$ Due to the convergence of the
approximations $F^m(u_{0\delta},\lambda)$ to the fixed point
$u_{\delta}=u_{\delta}(\lambda)$ of the operator
$F(\cdot,\lambda):C((-\infty,\delta],BC(\Omega,R^n))\rightarrow
C((-\infty,\delta],BC(\Omega,R^n))$ we get
$\|{u}_{\delta}(\lambda)-u_{0\delta}\|_{C((-\infty,\delta],BC(\Omega,R^n))}\leq
\varepsilon$ for each $\lambda\in U_0\bigcap U_1$ and
$\varepsilon\to0$ as $\lambda\to\lambda_0$.

Now,  addressing the model $(3.3)$ and Theorem 3.2, and taking
$\varphi_\lambda={u}_{\delta}(\lambda)$ and $a=\delta$, we prove
this theorem. \end{proof}

We note here that the remark, analogous to  Remark 3.3, is valid for
Theorem 3.4 as well.

\null

\begin{remark}If $\Omega$ is bounded, we can get the conclusion of   Theorem 3.4  replacing
3) and 4) by the following conditions:

\vspace{10pt}

 For any $\{\lambda_i\}\subset\Lambda$,
$\lambda_i\to\lambda_0$ it holds true that:

\vspace{10pt}

 3$^\ast$)  For any  $u\in R^n$ we have
$|f_{\lambda_i}(u)-f_{\lambda_0}(u)|\to 0$;

\vspace{10pt}

 4$^\ast$)  For all $x\in\Omega$,
$|\tau_{\lambda_i}(\cdot,x,\cdot)-\tau_{\lambda_0}(\cdot,x,\cdot
)|\to 0$ in measure on $(-\infty, b]\times\Omega$.
\end{remark}

\null

 Proof of the statement in Remark 3.7 is given in Appendix C.

\null

In case of spatially localized solutions to the $(1.7)$ and $(3.4)$,
we have the following remark to Theorem 3.4.

\null

\begin{remark}If in $(3.4)$ we add the condition, corresponding to localization
 in the spatial variable,

 \vspace{10pt}

$(\mathcal{A}_\lambda^\prime6)$   For each $\lambda\in \Lambda$ and
any $b\in R$, $\lim\limits_{|x|\to\infty}\big|
W_\lambda(t,s,x,y)\big|   = 0$ for all
 $(t,s,y)\in (-\infty,b]\times(-\infty,b]\times\Omega$,
 then, in order to get the conclusion of   Theorem 3.4  for
spatially localized solutions,  we need the following conditions
  instead of   2), 3), and 4), respectively:

\vspace{10pt}

 For any $\{\lambda_i\}\subset\Lambda$,
$\lambda_i\to\lambda_0$ it holds true that:

\vspace{10pt}

2$^\prime$)  For any $b\in R$, $r>0$ and each $t\in(-\infty,b]$,
$x\in \Omega, \ |x|\leq r$ it holds true that
$$\Big|\int\limits_{-\infty}^t\int\limits_{\Omega}\Big(W_{\lambda_i}(t,s,x,y)-
W_{\lambda_0}(t,s,x,y)\Big)dyds\Big|\to 0;$$   \vspace{10pt}

 3$^\prime$)  For any  $u\in R^n$ we have
$|f_{\lambda_i}(u)-f_{\lambda_0}(u)|\to 0$;

\vspace{10pt}

 4$^\prime$)  For all $x\in\Omega$,
$|\tau_{\lambda_i}(\cdot,x,\cdot)-\tau_{\lambda_0}(\cdot,x,\cdot
)|\to 0$ in measure on $(-\infty, b]\times\Omega$.
\end{remark}

\null

 Proof of the statement in Remark 3.8 is given in Appendix D.

\medskip

\section{ Conclusions and Outlook }

For the nonlinear Volterra integral equations $(1.7)$ and $(3.1)$,
which generalize the commonly used in the neural field theory models
$(1.1)$ -- $(1.6)$, we have defined the notions of local, global and
maximally extended solutions. We have obtained conditions which
guarantee existence of a unique global or maximally extended
solution and its continuous dependence on the equation parameters.
These results can also serve as a starting point  for the
development of numerical schemes for a broad class of neural field
models. A key word in this context is justification of such schemes.
We will emphasize that our results shed light on the problem of
structural stability in nonlocal field models in, e.g. systems
biology.

\null

 \verb"This is a draft of the paper containing the main results with the proofs."

\verb"Full-text version is available at "

http://math-res-pub.org/jadea/6/1/wellposedness-generalized-neural-field-equations-delay

\null

%

\medskip

\section*{Appendix A. \ Proof of The Statement in Remark 3.4}

We refer here to the proof of Theorem 3.2 and note that conditions
in
 Remark 3.4 imply that
$$|(S_{\tau_\lambda}^{\varphi_\lambda}
u_i)(\cdot,x,\cdot)-(S_{\tau_{\lambda_0}}^{\varphi_{\lambda_0}}
u)(\cdot,x,\cdot)|\to 0$$ uniformly on  $\big(([a,b]
\times\Omega)\setminus \Theta_{\lambda}\big)\times\ R^n$
($\mu(\Theta_{\lambda})\to0$), for each $x\in\Omega$, as
$\|u_i-u\|_{C([a,b],C_{0}(\Omega,R^n))}\to0$ and
$\lambda\to\lambda_0$.

Choose arbitrary $\varepsilon>0$. For the $b$ chosen in the proof of
Theorem 3.2 we find $$a^\ast=\min\limits_{t\in[a,b];\
(x,y)\in\Omega^2}(t-\tau_\lambda(t,x,y)).$$ Define the piecewise
constant functions $\overline{u}: [a, b]\times {R}^n \to {R}^n$ and
$\overline{\varphi}_{\lambda_0}: [a^\ast, a]\times {R}^n \to {R}^n$
as $\overline{u}(t,x)\in {R}^n$  for $t\in[a,b]$, $\xi\in [a^\ast,
a]$, $x\in\Omega$ such that
$$
\left\{%
\begin{array}{ll}
    |\overline{u}(t,x)-{u}(t,x)|\leq \varepsilon/2, \ \mbox{if} \ \ |\overline{u}(t,x)|>|{u}(t,x)|;\\
    |\overline{u}(t,x)-{u}(t,x)|< \varepsilon/2, \ \mbox{if} \ \ |\overline{u}(t,x)|<|{u}(t,x)|;\\
\end{array}%
\right.
$$
$$
\left\{%
\begin{array}{ll}
    |\overline{\varphi}_{\lambda_0}(\xi,x)-{\varphi}_{\lambda_0}(\xi,x)|\leq \varepsilon/2, \ \mbox{if} \ \ |\overline{\varphi}_{\lambda_0}(\xi,x)|>|{\varphi}_{\lambda_0}(\xi,x)|;\\
    |\overline{\varphi}_{\lambda_0}(\xi,x)-{\varphi}_{\lambda_0}(\xi,x)|< \varepsilon/2, \ \mbox{if} \ \ |\overline{\varphi}_{\lambda_0}(\xi,x)|<|{\varphi}_{\lambda_0}(\xi,x)|.\\
\end{array}%
\right.
$$
  We get the estimate
$$
\Big|f_\lambda\Big(\big(S_{\tau_\lambda}^{\varphi_\lambda} u_i
\big)(t,x,y)\Big)-f_{\lambda_0}\Big(\big(S_{\tau_{\lambda_0}}^{\varphi_{\lambda_0}}
u  \big)(t,x,y)\Big)\Big|\leq$$
$$\leq\Big|f_\lambda\Big(\big(S_{\tau_\lambda}^{\varphi_\lambda} u_i
\big)(t,x,y)\Big)-f_{\lambda}\Big(\big(S_{\tau_{\lambda_0}}^{\overline{\varphi}_{\lambda_0}}
\overline{u}  \big)(t,x,y)\Big)\Big|+
$$
$$
+\Big|f_{\lambda}\Big(\big(S_{\tau_{\lambda_0}}^{\overline{\varphi}_{\lambda_0}}
\overline{u}  \big)(t,x,y)\Big)
-f_{\lambda_0}\Big(\big(S_{\tau_{\lambda_0}}^{\overline{\varphi}_{\lambda_0}}
\overline{u}  \big)(t,x,y)\Big)\Big|+
$$
$$
+\Big|f_{\lambda_0}\Big(\big(S_{\tau_{\lambda_0}}^{\overline{\varphi}_{\lambda_0}}
\overline{u}
\big)(t,x,y)\Big)-f_{\lambda_0}\Big(\big(S_{\tau_{\lambda_0}}^{\varphi_{\lambda_0}}
u  \big)(t,x,y)\Big)\Big|.
$$
  Using the functions
$\overline{u}$ and $\overline{\varphi}_{\lambda_0}$, it is easy to
conclude  that the first and the third terms  on the right-hand side
of this inequality are less or equal to $2\varepsilon$ and
$\varepsilon$, respectively, on $\big(([a,b] \times\Omega)\setminus
\Theta_{\lambda}\big)\times\Omega$, where
$\mu(\Theta_{\lambda})\to0$ as  $\lambda\to\lambda_0$. In addition,
the condition 4$^\ast$) provide convergence to $0$ of the second
term on the right-hand side of the inequality as
$\lambda\to\lambda_0$.

Using the convergence obtained above, we get
$$
\max\limits_{t\in[a,b],\ x\in\Omega}\Big|\int\limits_a^t
 \int\limits_{\Omega}W_{\lambda}(t,s,x,y)f_\lambda\Big(\big(S^{\varphi_{\lambda}}_{\tau_{\lambda}}
u_i\big)(s,x,y)\Big)dyds- \ \ \ \ \ \ \ \ \ \ \ \ \ \ \ \ \ \ \ \ \
\
$$$$\ \ \ \ \ \ \ \ \ \ \ \ \
\ \ \ \ \ \ \ \-\int\limits_a^t
 \int\limits_{\Omega}W_{\lambda_0}(t,s,x,y)f_{\lambda_0}\Big(\big(S^{\varphi_{\lambda_0}}_{\tau_{\lambda_0}}
u\big)(s,x,y)\Big)dyds\Big|\leq
$$
$$
\max\limits_{t\in[a,b],\ x\in\Omega}\Big|\int\limits_a^t
 \int\limits_{\Omega}W_{\lambda}(t,s,x,y)f_\lambda\Big(\big(S^{\varphi_{\lambda}}_{\tau_{\lambda}}
u_i\big)(s,x,y)\Big)dyds- \ \ \ \ \ \ \ \ \ \ \ \ \ \ \ \ \ \ \ \ \
\
$$$$\ \ \ \ \ \ \ \ \ \ \ \ \
\ \ \ \ \ \ \ \-\int\limits_a^t
 \int\limits_{\Omega}W_{\lambda}(t,s,x,y)f_{\lambda_0}\Big(\big(S^{\varphi_{\lambda_0}}_{\tau_{\lambda_0}}
u\big)(s,x,y)\Big)dyds\Big|+
$$
$$
+\max\limits_{t\in[a,b],\ x\in\Omega}\Big|\int\limits_a^t
 \int\limits_{\Omega}W_{\lambda}(t,s,x,y)f_{\lambda_0}\Big(\big(S^{\varphi_{\lambda_0}}_{\tau_{\lambda_0}}
u\big)(s,x,y)\Big)dyds- \ \ \ \ \ \ \ \ \ \ \ \ \ \ \ \ \ \ \ \ \ \
$$$$\ \ \ \ \ \ \ \ \ \ \ \ \
\ \ \ \ \ \ \ \-\int\limits_a^t
 \int\limits_{\Omega}W_{\lambda_0}(t,s,x,y)f_{\lambda_0}\Big(\big(S^{\varphi_{\lambda_0}}_{\tau_{\lambda_0}}
u\big)(s,x,y)\Big)dyds\Big|.
$$
Taking into account the condition 3$^\ast$), we have the first term
on the right-hand side of this inequality going to $0$ as
$\lambda\to\lambda_0$. The second   term  on the right-hand side of
the inequality goes to $0$ by the virtue of the condition 2)  as
$\lambda\to\lambda_0$. Thus, the statement in Remark 3.4 is valid.

\medskip

\section*{Appendix B. \ Proof of The Statement in  Remark 3.5 }

 Conditions in Remark 3.5 imply the following changes in the proof of Theorem
 3:
$$|(S_{\tau_\lambda}^{\varphi_\lambda}
u_i)(\cdot,x,\cdot)-(S_{\tau_{\lambda_0}}^{\varphi_{\lambda_0}}
u)(\cdot,x,\cdot)|\to 0$$ uniformly on  $\big(([a,b]
\times\Omega)\setminus \Theta_{\lambda}\big)\times\ R^n$
($\mu(\Theta_{\lambda})\to0$), for each $x\in\Omega$, as
$\|u_i-u\|_{C([a,b],C_{0}(\Omega,R^n))}\to0$ and
$\lambda\to\lambda_0$.

Choose arbitrary $\varepsilon>0$. Define the piecewise constant
functions $\overline{u}: [a, b]\times {R}^n \to {R}^n$ and
$\overline{\varphi}_{\lambda_0}: (-\infty, a]\times {R}^n \to {R}^n$
as $\overline{u}(t,x)\in {R}^n$  for $t\in[a,b]$, $\xi\in (-\infty,
a]$, $x\in\Omega$ such that
$$
\left\{%
\begin{array}{ll}
    |\overline{u}(t,x)-{u}(t,x)|\leq \varepsilon/2, \ \mbox{if} \ \ |\overline{u}(t,x)|>|{u}(t,x)|;\\
    |\overline{u}(t,x)-{u}(t,x)|< \varepsilon/2, \ \mbox{if} \ \ |\overline{u}(t,x)|<|{u}(t,x)|;\\
\end{array}%
\right.
$$
$$
\left\{%
\begin{array}{ll}
    |\overline{\varphi}_{\lambda_0}(\xi,x)-{\varphi}_{\lambda_0}(\xi,x)|\leq \varepsilon/2, \ \mbox{if} \ \ |\overline{\varphi}_{\lambda_0}(\xi,x)|>|{\varphi}_{\lambda_0}(\xi,x)|;\\
    |\overline{\varphi}_{\lambda_0}(\xi,x)-{\varphi}_{\lambda_0}(\xi,x)|< \varepsilon/2, \ \mbox{if} \ \ |\overline{\varphi}_{\lambda_0}(\xi,x)|<|{\varphi}_{\lambda_0}(\xi,x)|.\\
\end{array}%
\right.
$$
  We get the estimate
$$
\Big|f_\lambda\Big(\big(S_{\tau_\lambda}^{\varphi_\lambda} u_i
\big)(t,x,y)\Big)-f_{\lambda_0}\Big(\big(S_{\tau_{\lambda_0}}^{\varphi_{\lambda_0}}
u  \big)(t,x,y)\Big)\Big|\leq$$
$$\leq\Big|f_\lambda\Big(\big(S_{\tau_\lambda}^{\varphi_\lambda} u_i
\big)(t,x,y)\Big)-f_{\lambda}\Big(\big(S_{\tau_{\lambda_0}}^{\overline{\varphi}_{\lambda_0}}
\overline{u}  \big)(t,x,y)\Big)\Big|+
$$
$$
+\Big|f_{\lambda}\Big(\big(S_{\tau_{\lambda_0}}^{\overline{\varphi}_{\lambda_0}}
\overline{u}  \big)(t,x,y)\Big)
-f_{\lambda_0}\Big(\big(S_{\tau_{\lambda_0}}^{\overline{\varphi}_{\lambda_0}}
\overline{u}  \big)(t,x,y)\Big)\Big|+
$$
$$
+\Big|f_{\lambda_0}\Big(\big(S_{\tau_{\lambda_0}}^{\overline{\varphi}_{\lambda_0}}
\overline{u}
\big)(t,x,y)\Big)-f_{\lambda_0}\Big(\big(S_{\tau_{\lambda_0}}^{\varphi_{\lambda_0}}
u  \big)(t,x,y)\Big)\Big|.
$$
  Using the functions
$\overline{u}$ and $\overline{\varphi}_{\lambda_0}$, it is easy to
conclude  that the first and the third terms  on the right-hand side
of this inequality are less or equal to $2\varepsilon$ and
$\varepsilon$, respectively, on $\big(([a,b] \times\Omega)\setminus
\Theta_{\lambda}\big)\times\Omega$, where
$\mu(\Theta_{\lambda})\to0$ as  $\lambda\to\lambda_0$. In addition
to that, the condition 4$^\prime$) provide convergence to $0$ of the
second term on the right-hand side of the inequality as
$\lambda\to\lambda_0$.

Using the convergence obtained above, $(A_\lambda^\prime3)$,
$(A_\lambda5)$ , and conditions 2$^\prime$) and 3$^\prime$), we get
$$
\max\limits_{t\in[a,b],\ x\in\Omega}\Big|\int\limits_a^t
 \int\limits_{\Omega}W_{\lambda}(t,s,x,y)f_\lambda\Big(\big(S^{\varphi_{\lambda}}_{\tau_{\lambda}}
u_i\big)(s,x,y)\Big)dyds- \ \ \ \ \ \ \ \ \ \ \ \ \ \ \ \ \ \ \ \ \
\
$$$$\ \ \ \ \ \ \ \ \ \ \ \ \
\ \ \ \ \ \ \ \-\int\limits_a^t
 \int\limits_{\Omega}W_{\lambda_0}(t,s,x,y)f_{\lambda_0}\Big(\big(S^{\varphi_{\lambda_0}}_{\tau_{\lambda_0}}
u\big)(s,x,y)\Big)dyds\Big|\leq
$$
$$
\max\limits_{t\in[a,b],\ x\in\Omega}\Big|\int\limits_a^t
 \int\limits_{\{x\in \Omega, |x|\leq
r^\prime\}}W_{\lambda}(t,s,x,y)f_\lambda\Big(\big(S^{\varphi_{\lambda}}_{\tau_{\lambda}}
u_i\big)(s,x,y)\Big)dyds- \ \ \ \ \ \ \ \ \ \ \ \ \ \ \ \ \ \ \ \ \
\
$$$$\ \ \ \ \ \ \ \ \ \ \ \ \
\ \ \ \ \ \ \ \-\int\limits_a^t
 \int\limits_{\{x\in \Omega, |x|\leq
r^\prime\}}W_{\lambda}(t,s,x,y)f_{\lambda_0}\Big(\big(S^{\varphi_{\lambda_0}}_{\tau_{\lambda_0}}
u\big)(s,x,y)\Big)dyds\Big|+
$$
$$
+\max\limits_{t\in[a,b],\ x\in\Omega}\Big|\int\limits_a^t
 \int\limits_{\{x\in \Omega, |x|\leq
r^\prime\}}W_{\lambda}(t,s,x,y)f_{\lambda_0}\Big(\big(S^{\varphi_{\lambda_0}}_{\tau_{\lambda_0}}
u\big)(s,x,y)\Big)dyds- \ \ \ \ \ \ \ \ \ \ \ \ \ \ \ \ \ \ \ \ \ \
$$$$\ \ \ \ \ \ \ \ \ \ \ \ \
\ \ \ \ \ \ \ \-\int\limits_a^t
 \int\limits_{\{x\in \Omega, |x|\leq
r^\prime\}}W_{\lambda_0}(t,s,x,y)f_{\lambda_0}\Big(\big(S^{\varphi_{\lambda_0}}_{\tau_{\lambda_0}}
u\big)(s,x,y)\Big)dyds\Big|+\epsilon_{r^\prime}(t,x).
$$
Here $\epsilon_{r^\prime}(t,x)\to0$ uniformly as $r^\prime\to
\infty$. Taking into account the condition 3$^\prime$), we have the
first term on the right-hand side of this inequality going to $0$ as
$\lambda\to\lambda_0$. The second   term  on the right-hand side of
the inequality goes to $0$ by the virtue of the condition
2$^\prime$)  as $\lambda\to\lambda_0$. Thus, the statement in Remark
3.5 is valid.

\medskip

\section*{Appendix C. \ Proof of The Statement in Remark 3.7 }

 The following changes in the proof of Theorem 3.4 stem from the conditions of Remark
3.7:
$$|u_i(t-{\tau_\lambda}(t,x,y),y)-u(t-{\tau_{\lambda_0}}(t,x,y),y)|\to 0$$
uniformly on  $\big(((-\infty,b] \times\Omega)\setminus
\Theta_{\lambda}\big)\times\ R^n$ ($\mu(\Theta_{\lambda})\to0$) for
each $x\in\Omega$, as
$\|u_i-u\|_{C((-\infty,b],BC(\Omega,R^n))}\to0$ and
$\lambda\to\lambda_0$.

Choose an arbitrary $\varepsilon>0$. Define the piecewise constant
function  $\overline{u}: (-\infty,b]\times {R}^n \to {R}^n$ as
$\overline{u}(t,x)\in {R}^n$  for $t\in(-\infty,b]$, $x\in\Omega$
such that
$$
\left\{%
\begin{array}{ll}
    |\overline{u}(t,x)-{u}(t,x)|\leq \varepsilon/2, \ \mbox{if} \ \ |\overline{u}(t,x)|>|{u}(t,x)|;\\
    |\overline{u}(t,x)-{u}(t,x)|< \varepsilon/2, \ \mbox{if} \ \ |\overline{u}(t,x)|<|{u}(t,x)|.\\
\end{array}%
\right.
$$
  Using the function introduced above, we get the estimate
$$
\Big|f_\lambda\Big(u_i(t-{\tau_\lambda}(t,x,y),y)\Big)-
f_{\lambda_0}\Big(u(t-{\tau_{\lambda_0}}(t,x,y),y)\Big)\Big|\leq$$
$$\leq\Big|f_\lambda\Big(u_i(t-{\tau_\lambda}(t,x,y),y)\Big)-
f_{\lambda}\Big(\overline{u}(t-{\tau_{\lambda_0}}(t,x,y),y)\Big)\Big|+
$$
$$
+\Big|f_{\lambda}\Big(\overline{u}(t-{\tau_{\lambda_0}}(t,x,y),y)\Big)
-f_{\lambda_0}\Big(\overline{u}(t-{\tau_{\lambda_0}}(t,x,y),y)\Big)\Big|+
$$
$$
+\Big|f_{\lambda_0}\Big(\overline{u}(t-{\tau_{\lambda_0}}(t,x,y),y)\Big)-f_{\lambda_0}\Big(u(t-{\tau_{\lambda_0}}(t,x,y),y)\Big)\Big|.
$$
  Here, the first and the third
terms  on the right-hand side of this inequality are less or equal
to $2\varepsilon$ and $\varepsilon$, respectively, on
$\big(((-\infty,b] \times\Omega)\setminus
\Theta_{\lambda}\big)\times\ R^n$, where $\mu(\Theta_{\lambda})\to0$
as  $\lambda\to\lambda_0$. In addition, the condition 4$^\ast$)
provide convergence to $0$ of the second term on the right-hand side
of the inequality as $\lambda\to\lambda_0$.

Using the convergence obtained above and $(\mathcal{A}_\lambda4)$,
we get
$$
\max\limits_{t\in(-\infty,b],\
x\in\Omega}\Big|\int\limits_{-\infty}^t
 \int\limits_{\Omega}W_{\lambda}(t,s,x,y)f_\lambda\Big(u_i(t-{\tau_\lambda}(t,x,y),y)\Big)dyds-
\ \ \ \ \ \ \ \ \ \ \ \ \ \ \ \ \ \ \ \ \ \
$$$$\ \ \ \ \ \ \ \ \ \ \ \ \
\ \ \ \ \ \ \ \-\int\limits_{-\infty}^t
 \int\limits_{\Omega}W_{\lambda_0}(t,s,x,y)
f_{\lambda_0}\Big(u(t-{\tau_{\lambda_0}}(t,x,y),y)\Big)dyds\Big|\leq
$$
$$
\max\limits_{t\in(-\infty,b],\
x\in\Omega}\Big|\int\limits_{-\infty}^t
 \int\limits_{\Omega}W_{\lambda}(t,s,x,y)f_\lambda\Big(u_i(t-{\tau_\lambda}(t,x,y),y)\Big)dyds-
\ \ \ \ \ \ \ \ \ \ \ \ \ \ \ \ \ \ \ \ \ \
$$$$\ \ \ \ \ \ \ \ \ \ \ \ \
\ \ \ \ \ \ \ \-\int\limits_{-\infty}^t
 \int\limits_{\Omega}W_{\lambda}(t,s,x,y)f_{\lambda_0}\Big(u(t-{\tau_{\lambda_0}}(t,x,y),y)\Big)dyds\Big|+
$$
$$
+\max\limits_{t\in(-\infty,b],\
x\in\Omega}\Big|\int\limits_{-\infty}^t
 \int\limits_{\Omega}W_{\lambda}(t,s,x,y)f_{\lambda_0}\Big(u(t-{\tau_{\lambda_0}}(t,x,y),y)\Big)dyds-
\ \ \ \ \ \ \ \ \ \ \ \ \ \ \ \ \ \ \ \ \ \
$$$$\ \ \ \ \ \ \ \ \ \ \ \ \
\ \ \ \ \ \ \ \-\int\limits_{-\infty}^t
 \int\limits_{\Omega}W_{\lambda_0}(t,s,x,y)f_{\lambda_0}\Big(u(t-{\tau_{\lambda_0}}(t,x,y),y)\Big)dyds\Big|.
$$
 Taking into account the condition 3$^\ast$), we have the
first term on the right-hand side of this inequality going to $0$ as
$\lambda\to\lambda_0$. The second term  on the right-hand side of
the inequality goes to $0$ by the virtue of the conditions 2) as
$\lambda\to\lambda_0$. Thus, the statement in Remark 3.7 is valid.

\medskip

\section*{Appendix D. \ Proof of The Statement in  Remark 3.8}

Referring to the proof of Theorem 3.4 we get the following changes
caused by conditions of Remark 3.8:
$$|u_i(t-{\tau_\lambda}(t,x,y),y)-u(t-{\tau_{\lambda_0}}(t,x,y),y)|\to 0$$
uniformly on  $\big(((-\infty,b] \times\Omega)\setminus
\Theta_{\lambda}\big)\times\ R^n$ ($\mu(\Theta_{\lambda})\to0$) for
each $x\in\Omega$, as
$\|u_i-u\|_{C((-\infty,b],BC(\Omega,R^n))}\to0$ and
$\lambda\to\lambda_0$.

Choose an arbitrary $\varepsilon>0$. Define the piecewise constant
function  $\overline{u}: (-\infty,b]\times {R}^n \to {R}^n$ as
$\overline{u}(t,x)\in {R}^n$  for $t\in(-\infty,b]$, $x\in\Omega$
such that
$$
\left\{%
\begin{array}{ll}
    |\overline{u}(t,x)-{u}(t,x)|\leq \varepsilon/2, \ \mbox{if} \ \ |\overline{u}(t,x)|>|{u}(t,x)|;\\
    |\overline{u}(t,x)-{u}(t,x)|< \varepsilon/2, \ \mbox{if} \ \ |\overline{u}(t,x)|<|{u}(t,x)|.\\
\end{array}%
\right.
$$
  Using this function, we get the estimate
$$
\Big|f_\lambda\Big(u_i(t-{\tau_\lambda}(t,x,y),y)\Big)-
f_{\lambda_0}\Big(u(t-{\tau_{\lambda_0}}(t,x,y),y)\Big)\Big|\leq$$
$$\leq\Big|f_\lambda\Big(u_i(t-{\tau_\lambda}(t,x,y),y)\Big)-
f_{\lambda}\Big(\overline{u}(t-{\tau_{\lambda_0}}(t,x,y),y)\Big)\Big|+
$$
$$
+\Big|f_{\lambda}\Big(\overline{u}(t-{\tau_{\lambda_0}}(t,x,y),y)\Big)
-f_{\lambda_0}\Big(\overline{u}(t-{\tau_{\lambda_0}}(t,x,y),y)\Big)\Big|+
$$
$$
+\Big|f_{\lambda_0}\Big(\overline{u}(t-{\tau_{\lambda_0}}(t,x,y),y)\Big)-f_{\lambda_0}\Big(u(t-{\tau_{\lambda_0}}(t,x,y),y)\Big)\Big|.
$$
  Using the function
$\overline{u}$, it is easy to conclude that the first and the third
terms  on the right-hand side of this inequality are less or equal
to $2\varepsilon$ and $\varepsilon$, respectively, on
$\big(((-\infty,b] \times\Omega)\setminus
\Theta_{\lambda}\big)\times\ R^n$, where $\mu(\Theta_{\lambda})\to0$
as  $\lambda\to\lambda_0$. In addition, the condition 4$^\prime$)
provide convergence to $0$ of the second term on the right-hand side
of the inequality as $\lambda\to\lambda_0$.

Using the convergence obtained above,
$(\mathcal{A}_\lambda^\prime3)$, $(\mathcal{A}_\lambda4)$, and
conditions 2$^\prime$) and 3$^\prime$), we get
$$
\max\limits_{t\in(-\infty,b],\
x\in\Omega}\Big|\int\limits_{-\infty}^t
 \int\limits_{\Omega}W_{\lambda}(t,s,x,y)f_\lambda\Big(u_i(t-{\tau_\lambda}(t,x,y),y)\Big)dyds-
\ \ \ \ \ \ \ \ \ \ \ \ \ \ \ \ \ \ \ \ \ \
$$$$\ \ \ \ \ \ \ \ \ \ \ \ \
\ \ \ \ \ \ \ \-\int\limits_{-\infty}^t
 \int\limits_{\Omega}W_{\lambda_0}(t,s,x,y)
f_{\lambda_0}\Big(u(t-{\tau_{\lambda_0}}(t,x,y),y)\Big)dyds\Big|\leq
$$
$$
\max\limits_{t\in(-\infty,b],\
x\in\Omega}\Big|\int\limits_{-\infty}^t
 \int\limits_{\{x\in \Omega, |x|\leq
r^\prime\}}W_{\lambda}(t,s,x,y)f_\lambda\Big(u_i(t-{\tau_\lambda}(t,x,y),y)\Big)dyds-
\ \ \ \ \ \ \ \ \ \ \ \ \ \ \ \ \ \ \ \ \ \
$$$$\ \ \ \ \ \ \ \ \ \ \ \ \
\ \ \ \ \ \ \ \-\int\limits_{-\infty}^t
 \int\limits_{\{x\in \Omega, |x|\leq
r^\prime\}}W_{\lambda}(t,s,x,y)f_{\lambda_0}\Big(u(t-{\tau_{\lambda_0}}(t,x,y),y)\Big)dyds\Big|+
$$
$$
+\max\limits_{t\in(-\infty,b],\
x\in\Omega}\Big|\int\limits_{-\infty}^t
 \int\limits_{\{x\in \Omega, |x|\leq
r^\prime\}}W_{\lambda}(t,s,x,y)f_{\lambda_0}\Big(u(t-{\tau_{\lambda_0}}(t,x,y),y)\Big)dyds-
\ \ \ \ \ \ \ \ \ \ \ \ \ \ \
$$
$$\ \ \ \ \ \
\ \ \ \ \ \ \ \-\int\limits_{-\infty}^t
 \int\limits_{\{x\in \Omega, |x|\leq
r^\prime\}}W_{\lambda_0}(t,s,x,y)f_{\lambda_0}\Big(u(t-{\tau_{\lambda_0}}(t,x,y),y)\Big)dyds\Big|+\epsilon_{r^\prime}(t,x).
$$
Here $\epsilon_{r^\prime}(t,x)\to0$ uniformly as $r^\prime\to
\infty$. Taking into account the condition 3$^\prime$), we have the
first term on the right-hand side of this inequality going to $0$ as
$\lambda\to\lambda_0$. The second term  on the right-hand side of
the inequality goes to $0$ by the virtue of the conditions
2$^\prime$)   as $\lambda\to\lambda_0$. Thus, the statement in
Remark 3.8 is valid.

\label{lastpage-01}
\end{document}